\documentclass[twoside,leqno,twocolumn]{article}

\usepackage[letterpaper]{geometry}

\usepackage{ltexpprt}

\usepackage{hyperref}

\usepackage{algpseudocode} 

\usepackage{graphicx}%
\usepackage{multirow}%
\usepackage{amsmath,amssymb,amsfonts}%
\usepackage{mathrsfs}%
\usepackage[title]{appendix}%
\usepackage{xcolor}%
\usepackage{textcomp}%
\usepackage{manyfoot}%
\usepackage{booktabs}%
\usepackage{algorithmicx}%
\usepackage{stfloats}
\usepackage{comment}

\def\ra{$r$-adaptivity}

\def\bx{{\boldsymbol{x}}}   

\def\lim{{\xi}}    
\def\bbx{\bar{\bx}}

\def\ua{\underline{\alpha}}

\definecolor{code}{rgb}{0.7, 0, 0.4}

\newcommand{\rev}[1]{\textcolor{black}{#1}}

\begin{document}

\newcommand\anonauth[2]{{#1}} 
\newcommand\anon[2]{{#2}} 
\newcommand\relatedversion{}
\renewcommand\relatedversion{\thanks{The full version of the paper can be accessed at \anon{\protect\url{redacted-url}}{\protect\url{https://arxiv.org/abs/1902.09310}}}} 

\title{High-Order Mesh $r$-Adaptivity with Tangential Relaxation and Guaranteed Mesh Validity}
{ \author{
Ketan Mittal\thanks{Lawrence Livermore National Laboratory, Livermore, CA, U.S.A. \{mittal3,dobrev1,kolev1,tomov2\}@llnl.gov.} \and \hspace{-10mm} \thanks{Corresponding Author.} \and
Veselin Dobrev\footnotemark[1] \and
Tzanio Kolev\footnotemark[1] \and
Vladimir Tomov\footnotemark[1]
}}

\date{}

\maketitle

\begin{abstract} \small\baselineskip=9pt
High-order meshes are crucial for achieving optimal convergence rates in curvilinear domains, preserving symmetry, and aligning with key flow features in moving mesh simulations \cite{Dobrev2012}, but their quality is challenging to control.
In prior work, we have developed techniques based on Target-Matrix Optimization Paradigm (TMOP) to adapt a given high-order mesh to the geometry and solution of the partial differential equation (PDE) \cite{TMOP2019SISC,TMOP2020CAF}.
Here, we extend this framework to address two key gaps in the literature for high-order mesh $r$-adaptivity.
First, we introduce tangential relaxation on curved surfaces using solely the discrete mesh representation, eliminating the need for access to underlying geometry (e.g., CAD model).
Second, we ensure a continuously positive Jacobian determinant throughout the domain. This determinant positivity is essential for using the high-order mesh resulting from $r$-adaptivity with arbitrary quadrature schemes in simulations.
The proposed approach is demonstrated to be robust using a variety of numerical experiments.
\end{abstract}

\section{Introduction} \label{sec:intro}

High-order (curved) finite-element and spectral-element meshes have become essential tools in the numerical solution of partial differential equations, particularly in applications demanding high accuracy per degree of freedom (e.g., computational fluid dynamics, structural mechanics, and electromagnetics). By representing the geometry with polynomials of degree $p>1$, high-order meshes ensure that a curved geometry of interest is modeled accurately.

To ensure an accurate numerical solution of the PDE of interest, mesh elements must have good quality in terms of shape, size, and orientation. Poor element shapes can lead to ill-conditioned algebraic systems \cite{mittal2019mesh} and, even worse, tangled (invalid) elements can cause solvers to fail altogether. Thus, controlling mesh quality and guaranteeing mesh validity is essential for leveraging the theoretical advantages of high-order meshes.

Despite significant progress in high-order mesh optimization \cite{alauzet2016decade,aparicio2023combining,frey2005anisotropic,fidkowski2011review,loseille2010fully,marcon2019variational,odier2021mesh,yano2012optimization,zhang2018curvilinear}, two key challenges remain: (i) guaranteeing mesh validity, i.e., ensuring that the Jacobian determinant is continuously positive everywhere in the domain, and (ii) geometrically accurate tangential relaxation of nodes on curved boundaries without requiring explicit access to the underlying geometry (e.g., CAD model).

Existing $r$-adaptivity techniques determine mesh validity by sampling the Jacobian determinant at a finite set of quadrature points \cite{dobrev2019target,sanjaya2020comparison,aparicio2024defining} or employ Bernstein bases to use the corresponding coefficients as bounds \cite{coppeans2024anisotropic,luo2002p,johnen2013geometrical,rochery2025metris}. Sampling can miss local inversions between quadrature points, while Bernstein bounds tend to be overly conservative \cite{dzanic2025method}. We guarantee mesh validity by computing provable bounds on the Jacobian determinant using recently developed techniques \cite{mittal2025general,dzanic2025method},
which circumvents the shortcomings of existing methods.
A lower bound on the Jacobian determinant is also required for mesh untangling based on shifted-barrier metrics \cite{knupp2022worst}. With these methods, mesh untangling can fail if the barrier is not low enough (due to undersampling), or be slow and inefficient if the barrier is too conservative (due to overestimation such as from Bernstein coefficients). We use provable lower bounds on the Jacobian determinant with shifted-barrier metrics, thus increasing the robustness of mesh untangling process.

To maximize mesh improvement during $r$-adaptivity, boundary nodes should be allowed to move tangentially along the surface without compromising geometric accuracy. Existing techniques for tangential relaxation typically leverage CAD parameterizations \cite{kirilov2024high, ruiz2016high,shivanna2010analytical,garimella2004triangular}, but these techniques require development of specialized methods in existing finite element- and spectral element frameworks to interact with the CAD model.
Other techniques for tangential relaxation include use of implicitly defined surfaces through distance function \cite{mittal2019mesh} and level-sets \cite{TMOP2021IMR,TMOP2023CAD}, but these approaches do not guarantee that the target geometry is accurately recovered. We enable tangential relaxation on curved surfaces by allowing boundary nodes to first move freely for mesh quality improvement, immediately followed by closest point projection on the original high-order mesh surface after each $r$-adaptivity iteration.
This approach thus does not require explicit access to underlying geometry, such as a CAD model, and guarantees that the optimized mesh surface conforms with the original high-order mesh.

\begin{figure*}[tb]
\begin{center}
$\begin{array}{cc}
    \includegraphics[height=0.4\textwidth]{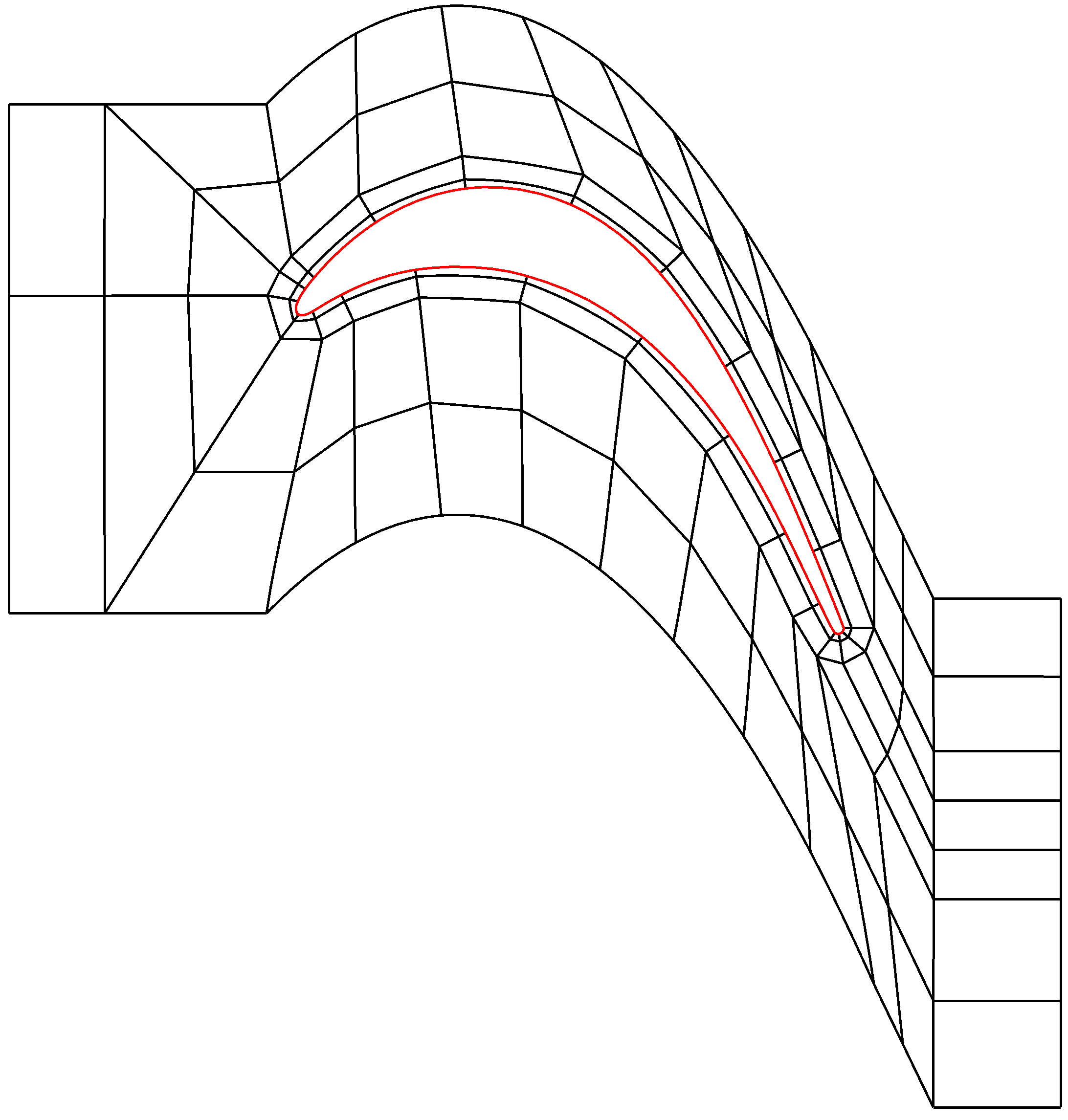} &
    \includegraphics[height=0.4\textwidth]{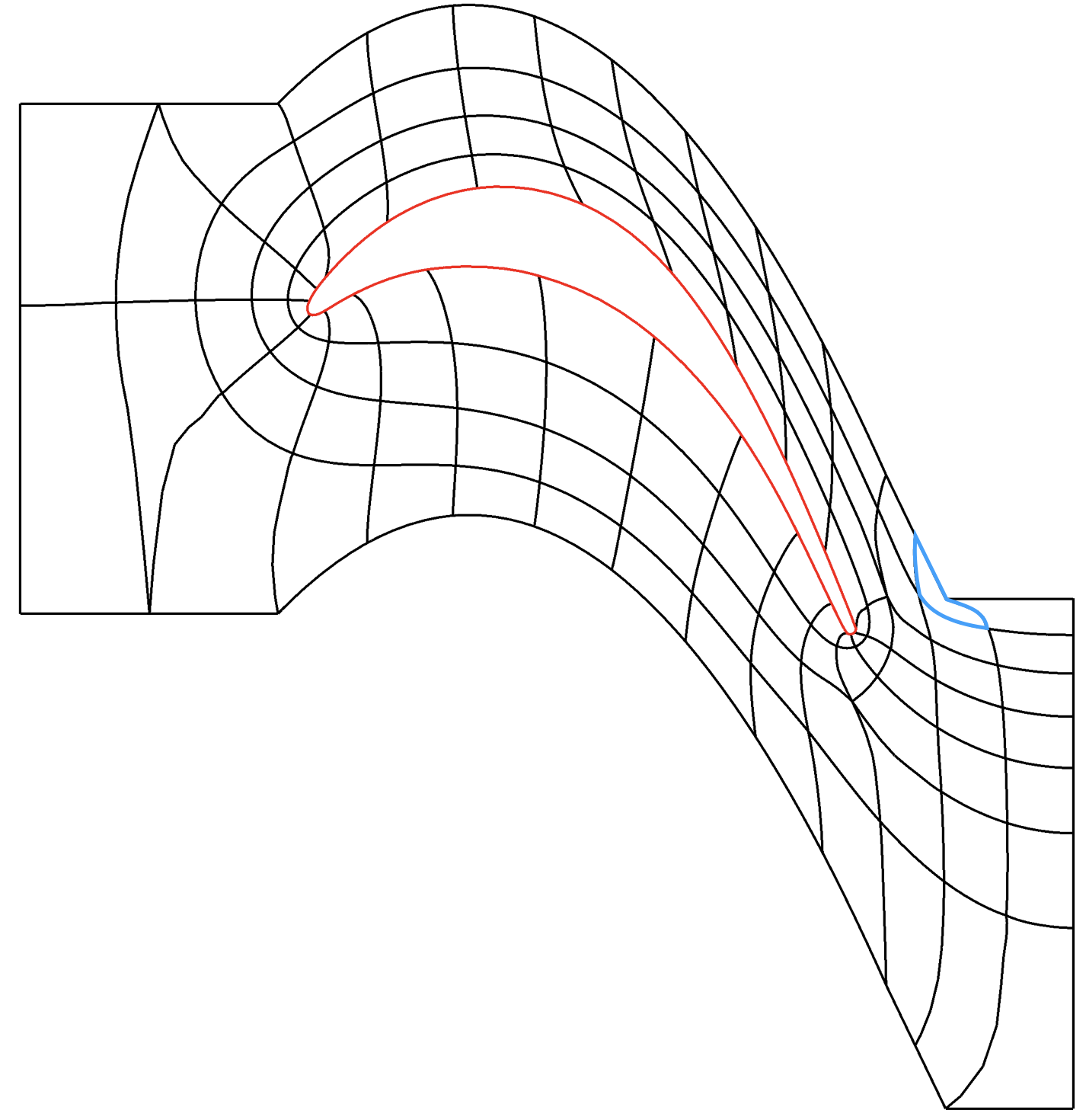} \\
    \textrm{(a) Initial mesh} & \textrm{(b) Optimized mesh}
\end{array}$
\end{center}
\caption{Blade mesh example illustrating (a) initial high-order mesh and (b) mesh after optimization.}
\label{fig:blade_example}
\end{figure*}

Figure \ref{fig:blade_example} illustrates the shortcomings of existing techniques in the context of $r$-adaptivity for a fourth-order mesh used to simulate flow around a turbine blade.
Since the boundary nodes are held fixed along the surface of the blade, the optimized mesh has skewed elements next to the upper surface.
Additionally, since the determinant of the Jacobian is only checked at a set of quadrature points, one of the elements becomes tangled (highlighted in blue) even though the minimum sampled Jacobian is positive ($2.8 \times 10^{-5}$).  Note that if Bernstein bases are used to represent the 7th-order function for the Jacobian determinant of the problematic element, the estimated lower bound is -0.0302, which is an order of magnitude below the actual minimum of -0.00104.

The remainder of the paper is organized as follows. Section \ref{sec:prelim} describes our current approach for $r$-adaptivity. Section \ref{sec:method} presents the proposed approach for guaranteeing mesh validity and tangential relaxation on curved surfaces. The effectiveness of these methods is demonstrated through various numerical experiments in Section \ref{sec:results}.
Summary and directions for future work are given in Section \ref{sec:conc}.

\section{Preliminaries} \label{sec:prelim}

\subsection{High-order mesh and function representation} \label{subsec:mesh}

In our finite-element-based framework, the domain $\Omega \subset \mathbb{R}^d$
is discretized as a union of curved mesh elements, $\Omega^e$, $e=1\dots N_E$,
each of order $p$.  To obtain a
discrete representation of these elements, we select a set of scalar basis
functions $\{ w_i \}_{i=1}^{N_p}$, on the $N_p$ nodes of the reference element $\bar{\Omega}^e$.
In the case of tensor-product elements (quadrilaterals in 2D and hexahedra in 3D),
$N_p = (p+1)^d$, and the basis spans the space
of all polynomials of degree at most $p$ in each variable.
These $p$th-order basis functions are typically chosen to be Lagrange
interpolants at the Gauss-Lobatto nodes of the reference element.
The position of an
element $\Omega^e$  in the mesh $\mathcal{M}$ is fully described by a matrix
$\bx_e$ of size $d \times N_p$ whose columns represent the coordinates
of the element {\em degrees of freedom} (DOFs).
Given $\bx_e$, we introduce the map between
the reference and physical element, $\Phi_e:\bar{\Omega}^e \to \mathbb{R}^d$:
\begin{equation}
\label{eq:x}
\vspace{-5mm}
\bx(\bbx) =
   \Phi_e(\bbx) \equiv
   \sum_{i=1}^{N_p} \mathbf{x}_{e,i} w_i(\bar{\bx}),
   \,\, \bbx \in \bar{\Omega}^e, ~ \bx \in \Omega^e,
\end{equation}
where $\mathbf{x}_{e,i}$ denotes the $i$-th column of $\bx_e$, i.e., the $i$-th node of element $\Omega^e$.

Throughout the manuscript, $\bx$ will denote the position function defined
by \eqref{eq:x}, and
$\bx_e$ will denote the element-wise vector/matrix of nodal locations for element $\Omega^e$.

Similar to \eqref{eq:x}, any high-order function $u$ is defined as
\begin{eqnarray}
\label{eq:u}
u(\bbx) = \sum_{j=1}^{N_u} u_j \phi_j(\bbx), \quad \bbx \in \bar{\Omega}^e
\end{eqnarray}
\rev{where $\{u_j\}_{j=1}^{N_u}$ are the coefficients for the bases $\{\phi_j\}_{j=1}^{N_u}$ and $N_u$ depends on the polynomial order used to represent the
function $u(\bbx)$.}

\subsection{TMOP for Mesh Quality Improvement via \ra}\label{subsec:tmop}

For a given element with nodal coordinates $\bx_e$, the Jacobian
of the mapping $\Phi_e$ at any reference point $\bbx \in \bar{\Omega}^e$ is
\begin{equation}
\label{eq:A}
  A_{ab}(\bbx) = \frac{\partial x_{a}(\bbx)}{\partial \bar{x}_b} =
    \sum_{i=1}^{N_p} x_{i,a} \frac{\partial \bar{w}_i(\bbx)}{\partial \bar{x}_b}, \quad a,b = 1 \dots d
\end{equation}
where $x_a$, represents the $a$th component of $\bx$ \eqref{eq:x}, and $x_{i,a}$ represents the $a$th component of $\bx_{e,i}$, i.e., the $i$th DOF in element $\Omega^e$.
The Jacobian matrix $A$ represents the local deformation of the physical element $\Omega^e$ with respect to the
reference element $\bar{\Omega}^e$ at the reference point $\bbx$. This matrix plays an important
role in FEM as it is used to determine mesh validity ($\alpha \equiv \det(A)$ must be greater than 0 at every point in the mesh),
and also used to compute derivatives and integrals. The Jacobian matrix can ultimately impact the accuracy
and computational cost of the solution \cite{mittal2019mesh}.
This lends to the central idea of TMOP of optimizing the mesh to control the
local Jacobian $A_{d \times d}$ in the mesh.

\subsubsection*{Target matrix} The first step for mesh optimization with TMOP is to specify a target transformation matrix $W_{d \times d}$, analogous to $A_{d \times d}$, for each point in the mesh. Target construction is typically guided by the fact that any Jacobian matrix can be written as a composition of four geometric components \cite{knupp2022geometric}, namely volume, rotation, skewness, and aspect ratio:
\begin{equation}
\label{eq:W}
W_{d\times d} = \underbrace{\zeta}_{\text{[volume]}} \circ
    \underbrace{R_{d\times d}}_{\text{[rotation]}} \circ
    \underbrace{Q_{d\times d}}_{\text{[skewness]}} \circ
    \underbrace{D_{d\times d}}_{\text{[aspect ratio]}}.
\end{equation}
In practice, the user may specify $W$ as a combination of any of the four fundamental components that they are interested in optimizing the mesh for. For the purposes of this paper, we are mainly concerned with ensuring good element \emph{shape} (skewness and aspect ratio), so we set $W$ to be that of an ideal element, i.e.,
square for quad elements, cube for hex elements, and
equilateral simplex for triangles and tetrahedrons.
Additionally, the size component $\zeta$ is held constant throughout the domain and is typically set based on the average volume of the elements in the mesh.
Advanced techniques on how spatially varying $W$ can be constructed for optimizing different geometric parameters, or even for automatically adapting the mesh to the solution of the PDE
are given in \cite{TMOP2020CAF,TMOP2021EWC,knupp2019target}.

\begin{figure}[tb!]
\centerline{\includegraphics[width=0.3\textwidth]{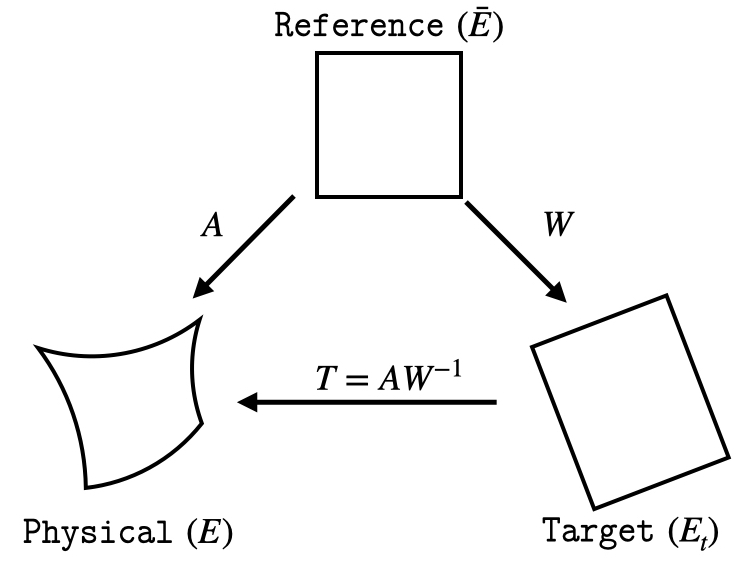}}
\vspace{-1mm}
\caption{Schematic representation of the major TMOP matrices.}
\label{fig:tmop}
\vspace{-3mm}
\end{figure}

\subsubsection*{Mesh quality metric} The next key component in the TMOP-based approach is a mesh quality
metric that measures the deviation between the current Jacobian transformation $A$
and the target transformation $W$.
The mesh quality metric $\mu(T)$, $T = A W^{-1}$ in Figure \ref{fig:tmop},
compares $A$ and $W$ in terms of some of the geometric parameters.
For example, $\mu_{2,s}=\frac{\mid T \mid^2}{2\tau}-1$ is a \emph{shape} metric\footnote{The metric subscript follows the numbering in \cite{Knupp2020, knupp2022geometric}.}
that depends on the skewness and aspect ratio components, but
is invariant to orientation/rotation and volume. Here, $|T|$ and $\tau$ are the Frobenius norm and determinant of $T$, respectively. Similarly,
$\mu_{77,v}=\frac{1}{2}(\tau-\tau^{-1})^2$ is a
\emph{size/volume} metric that depends only on the volume of the element.
Individual metrics can also be combined using a convex combination. For example, $\mu_{80,vs} =\gamma \mu_{2,s} + (1-\gamma) \mu_{77,v}$,
$0 \leq \gamma \leq 1$.
Occasionally, we also use metrics that directly depend on $A$ and $W$ (instead of via $T$) and we denote them as $\nu(A,W)$.

\subsubsection*{TMOP objective function} Using the mesh quality metric, the mesh optimization problem is minimizing the
global objective:
\begin{equation}
\label{eq:F_full}
    F(\bx) = \sum_{\Omega^e \in \mathcal{M}} \int_{\Omega^{e_t}} \mu(T(\bx)) d\bx_t,
\end{equation}
where $F$ is a sum of the TMOP objective function for each element in the mesh,
and $\Omega^{e_t}$ is the target element corresponding to the element $\Omega^e$.
Minimizing \eqref{eq:F_full} results in node movement such that
the local Jacobian transformation $A$ resembles the target transformation $W$ as close as possible at each point, in terms of the geometric parameters enforced by the mesh quality metric.

\subsubsection*{Solver} In our framework, we minimize $F$ iteratively using the Newton's method as
\begin{eqnarray}
\label{eq:r_adaptivity_solve}
\bx_{k+1} = \bx_k - \gamma \mathcal{H}^{-1}(\bx_k)\mathcal{J}(\bx_k).
\end{eqnarray}
Here, $\bx_k$ refers to the nodal positions at the $k$-th Newton iteration during \ra, $\gamma$ is a scalar determined by a line-search procedure, and
$\mathcal{H}(\bx_k)$ and $\mathcal{J}(\bx_k)$ are the
Hessian ($\partial^2 F(\bx_k) / \partial \bx_j\partial \bx_i$)
and the gradient ($\partial F(\bx_k) / \partial \bx_i$),
respectively, associated with the TMOP objective function.
The line-search procedure requires that $\gamma$ is chosen such that
\begin{subequations}\label{eq:line_search}
    \begin{align}
    \label{eq:line_search_1}
    F(\bx_{k+1}) < 1.2 F(\bx_k),\\
        \label{eq:line_search_2}
    |\mathcal{J}(\bx_{k+1})| < 1.2 |\mathcal{J}(\bx_k)|, \\
        \label{eq:line_search_3}
    \det\left( A\left( \bx_{k+1}(\bbx_e) \right) \right) > 0 \quad \forall \bbx_e \in \mathcal{Q}, \ \forall \Omega^e \in \mathcal{M}.
    \end{align}
\end{subequations}
\rev{The relaxed line-search conditions \eqref{eq:line_search_1}--\eqref{eq:line_search_2} have been tuned empirically to allow the solver to escape shallow local minima, which are common in mesh optimization landscapes.}
The last constraint in \eqref{eq:line_search_3} ensures that the mesh stays valid at the selected set of quadrature points $\mathcal{Q}$ in all elements at each iteration. \rev{We have found that these line-search conditions are sufficient for the robustness of the mesh optimizer for
our applications of interest \cite{TMOP2019SISC,TMOP2020CAF,TMOP2021EWC}.}
Comprehensive advances in solvers, preconditioners, and line-search procedures
for the same class of problems are presented in \cite{Roca2022}.

Finally, boundary nodes are allowed to move tangentially only on planar surfaces aligned with the Cartesian axes.
Newton iterations \eqref{eq:r_adaptivity_solve} are then
done until the line-search criterion cannot be satisfied or the relative $l_2$ norm of the gradient of the objective function
with respect to the current and original mesh nodal positions is below a certain
tolerance $\varepsilon$, i.e., $|\mathcal{J}(\bx)|/|\mathcal{J}(\bx_0)| \leq \varepsilon$. We set $\varepsilon = 10^{-10}$ for the results presented in the current work.

Figure \ref{fig:blade_example} shows the blade mesh optimized using this approach. The target matrix is set to identity ($W=I$) and a shape metric ($\mu_{2,s}$) is used. The resulting optimized mesh has elements closer to unity aspect ratio and skewness closer to $\pi/2$ radians in comparison to the original mesh, as prescribed by the target $W$.

\subsubsection*{Mesh untangling} The current $r$-adaptivity approach ensures that an initially valid mesh remains valid at quadrature points throughout the optimization process. When the initial mesh is invalid, a shifted-barrier metric is used \cite{knupp2022worst} in \eqref{eq:F_full}
\begin{eqnarray}
    \label{eq:shifted_barrier_qp}
    \mu(T) = \frac{\tilde{\mu}(T)}{2\big(\tau - \tau_b\big)}, \\
    \tau_{b} = \begin{cases}
    \beta \cdot \tau_{\tt qp,min}-\epsilon & \text{if } \tau_{\tt qp,min} \leq 0 \\
    0 & \text{otherwise}
    \end{cases}
\end{eqnarray}
where $\tilde{\mu}(T)$ is a \emph{non-barrier} metric (e.g., $\tilde{\mu}_4(T) = ||T||_F^2 - 2 \tau$), $\tau$ is the determinant of $T$, and $\tau_b$ is the barrier. The idea with shifted-barrier metrics is to have the barrier just below the minimum determinant ($\tau_{\tt min}$) across the mesh, such that $\mu(T)\rightarrow \infty$ as $\tau \rightarrow \tau_b$. Thus minimizing $\mu(T)$ results in mesh movement that moves $\tau$ away from the barrier, thus increasing $\tau_{\tt min}$ and eventually untangling the mesh.

Since $\tau_b$ is typically based on $\tau_{\tt qp,min}$, the minimum Jacobian determinant sampled across all the quadrature points ($\det(T(\bx_k))$), a scaling factor $\beta$ (typically set to 1.5) and offset $\epsilon$ (typically $10^{-2}$) are used to ensure that $\tau_b$ is below $\tau_{\tt min}$.
Note that the barrier is held at 0 once the mesh untangles. The condition $\tau_{\tt qp,min} > 0$ is equivalent to $\alpha_{\tt qp,min} > 0$ because the size component of the target matrix ($\omega \equiv \det(W)$) is positive throughout the domain.

\section{Methodology} \label{sec:method}

\subsection{Guaranteeing mesh validity using provable bounds on the Jacobian determinant} \label{subsec:bounds}

A general technique for bounding high-order finite element functions that use nodal or modal bases is described in \cite{mittal2025general} and \cite{dzanic2025method}.
The key idea therein is to construct piecewise linear bounds for each basis function in a given set of bases, and use them to bound any polynomial function expressed as a weighted sum of these basis functions (e.g., \eqref{eq:u}).

The piecewise linear bounds for each basis function ($\phi_i$, $i = 1\dots N$)  are determined by solving a constrained optimization problem that minimizes the $L_2$ norm of the difference between the basis function and the bounding function defined using a set of control nodes, $\eta_j \in \bar{\Omega}, j=1\dots M$, in the reference element. In the context of the upper bound, using $q_{+}^{ij}$ to denote the bounding values for the $i$th basis function at $\eta_j$, and $U^i_{\boldsymbol{\eta}, \mathbf{q}^{+}}$ to represent the vector of tuples denoting bounding values at the control nodes, the optimization problem is formulated as:
\begin{equation}
    \mathbf{q^{+}}^* = \underset{\mathbf{q^{+}}}{\arg \min}\ f(\mathbf{q^{+}}) \quad s.t.\quad  g(\mathbf{q^{+}}, \bbx) \geq 0 ,
\end{equation}
where
\begin{subequations}\label{eq:obj1}
    \begin{align}
    f(\mathbf{q^{+}}) &= \left \| U^i_{\boldsymbol{\eta}, \mathbf{q^{+}}}(\bbx) - \phi_i(\bbx)\right \|_{2, \Omega},\\
    g(\mathbf{q^{+}}, \bbx) &= U^i_{\boldsymbol{\eta}, \mathbf{q^{+}}}(\bbx) - \phi_i(\bbx).
    \end{align}
\end{subequations}
A similar optimization problem is solved to obtain lower bounding values, $q_{-}^{ij}$ for $L^i_{\boldsymbol{\eta}, \mathbf{q}^{-}}$.

One can then directly compute linear bounding functions for an arbitrary $u(\bbx)$ as
\begin{subequations}\label{eq:1D_bounds}
    \begin{align}
    \underline{u}(\eta_j) &= \sum_{i=1}^N \min(u_i q_{ij}^-,u_i q_{ij}^+),\\
    \overline{u}(\eta_j) &= \sum_{i=1}^N \max(u_i q_{ij}^-,u_i q_{ij}^+),
    \end{align}
\end{subequations}
such that $\underline{u} \leq u \leq \overline{u}$.
On each element $\Omega$, the lower and upper bound of the function $u(\bbx)$ is simply the pointwise minimum and maximum of $\underline{u}$ and $\overline{u}$ at all the control nodes, respectively. The effectiveness of this approach is further improved by optimizing for the control point locations, $\eta_j$, and removing a linear fit from the original function before bounding it. The reader is referred to \cite{mittal2025general,dzanic2025method} for a detailed description.

\begin{figure}[tb!]
\begin{center}
$\begin{array}{cc}
    \includegraphics[width=0.23\textwidth]{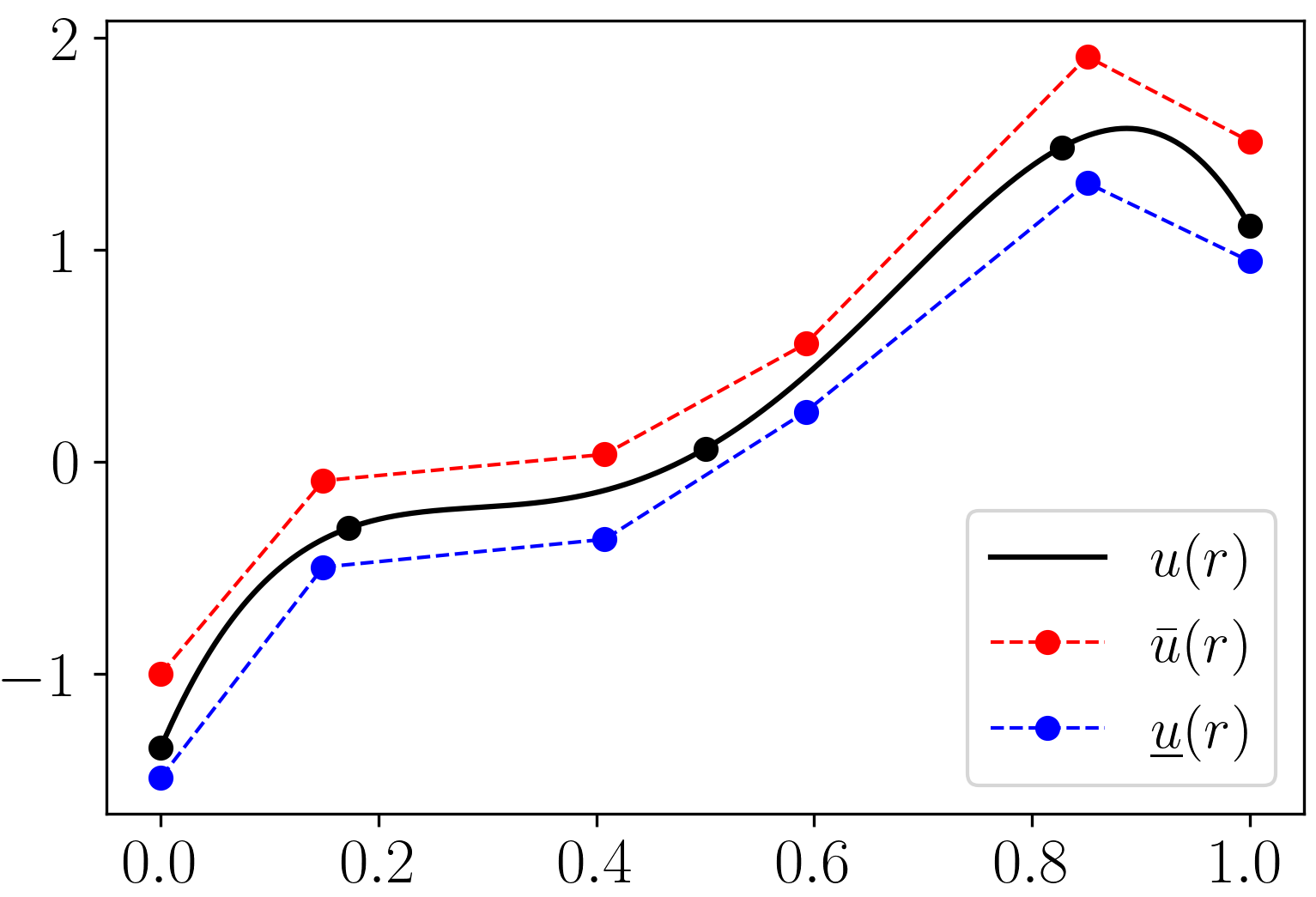} &
    \includegraphics[width=0.23\textwidth]{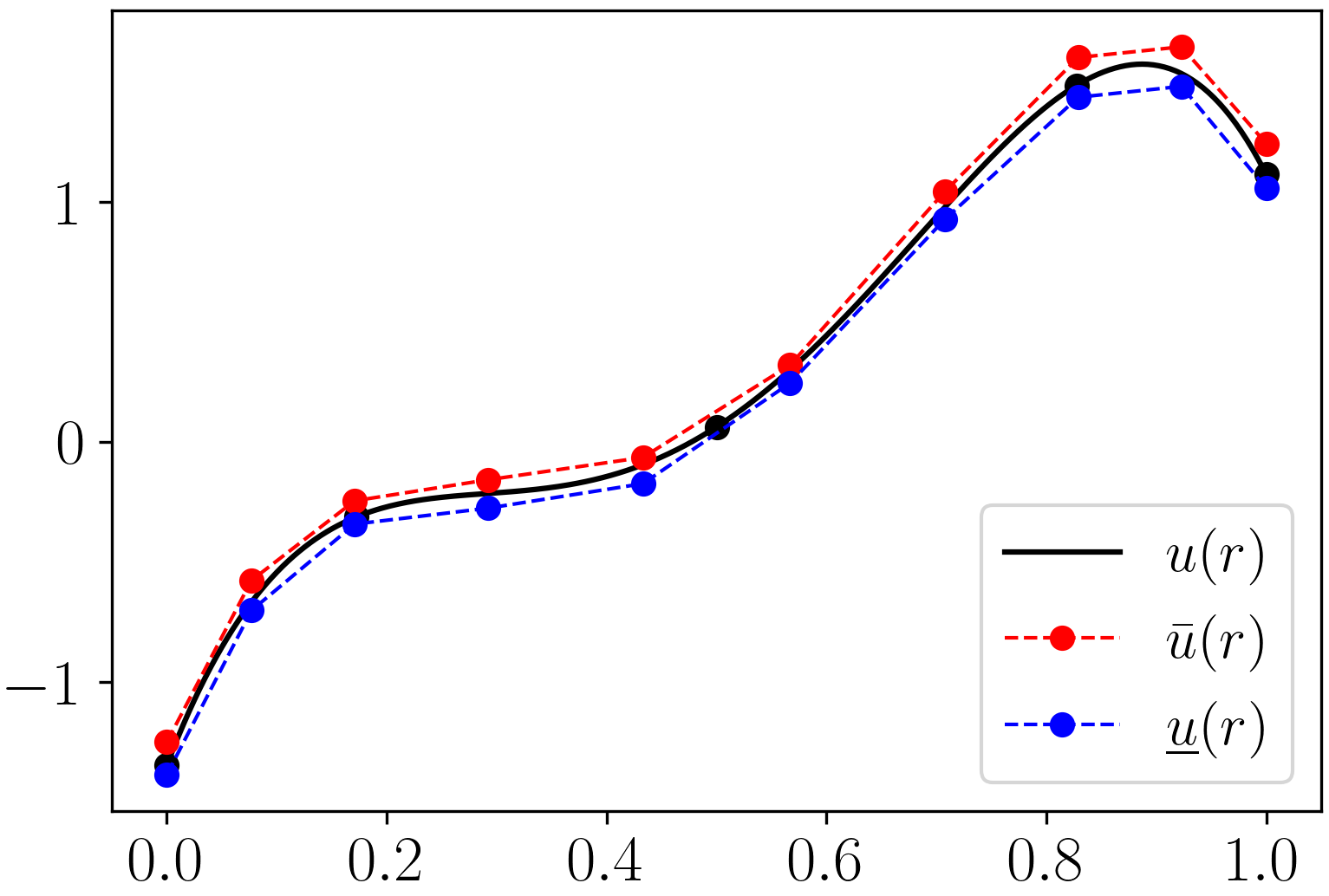} \\
\end{array}$
\end{center}
\vspace{-1mm}
\caption{Bounding a 4th-order function using the approach in \cite{mittal2025general,dzanic2025method} with (left) 6 control points and (right) 10 control points.}
\label{fig:bound_1D}
\vspace{-3mm}
\end{figure}

Figure \ref{fig:bound_1D} illustrates an example of this method to bound a 4th-order function $u(\bbx)$ with nodal coefficients $u = [-1.346, -0.311, 0.063, 1.485, 1.114]$ and Gauss-Lobatto-Legendre nodal bases.
Figure \ref{fig:bound_1D}(left) and (right) show the effectiveness of the proposed approach using different number of control points; more control points result in tighter bounds, as expected.

\subsubsection*{Mesh validity through bounds on the Jacobian determinant}

When the mesh nodal position function is of polynomial degree $p$, the determinant of its Jacobian is a polynomial of order $d\cdot p-1$ for tensor product elements and $d(p - 1)$ for simplices \cite{johnen2013geometrical}.
The preceding bounding procedure can thus be directly applied to compute $\underline{\alpha}$, the lower bound on the Jacobian determinant in the mesh.
We then guarantee that an initially valid mesh stays valid during $r$-adaptivity by accepting the Newton update \eqref{eq:r_adaptivity_solve} only if $\underline{\alpha}(\bx_{k+1}) > 0$.
Compactness of the bounds is ensured by using a sufficiently high number of control points or a recursive subdivision strategy.

\begin{figure}[tb!]
\begin{center}
$\begin{array}{ccc}
    \includegraphics[height=0.17\textwidth]{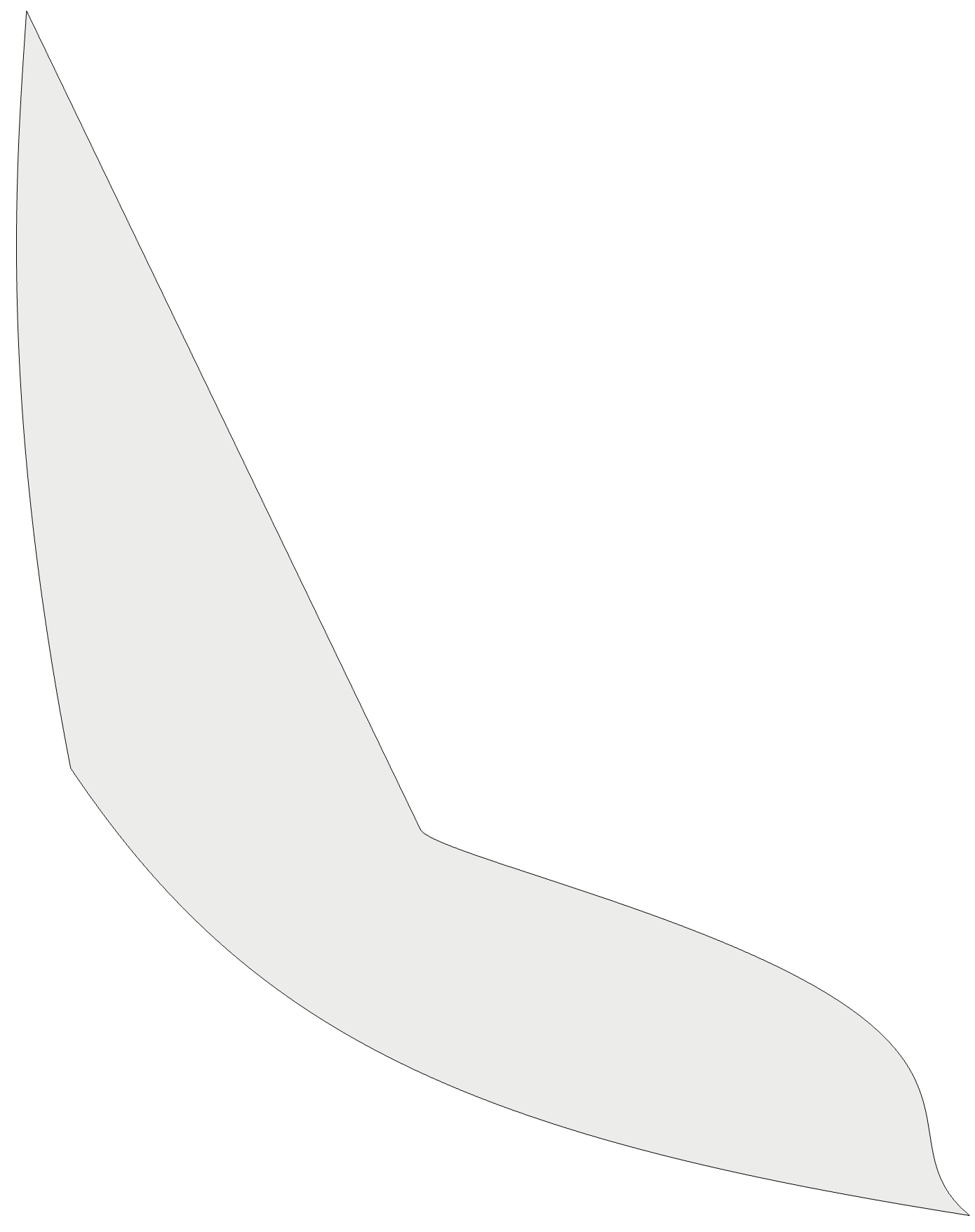} &
    \includegraphics[height=0.17\textwidth]{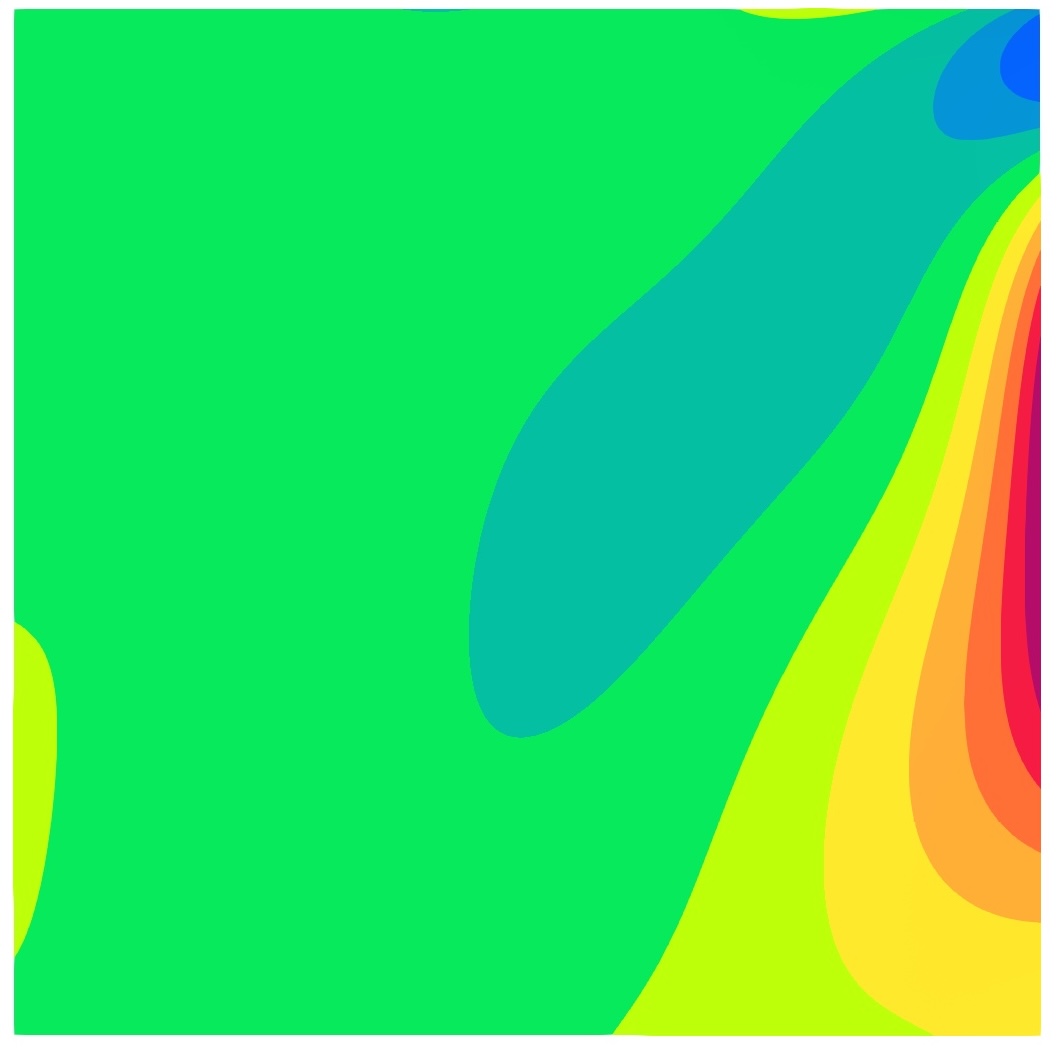} &
    \includegraphics[height=0.17\textwidth]{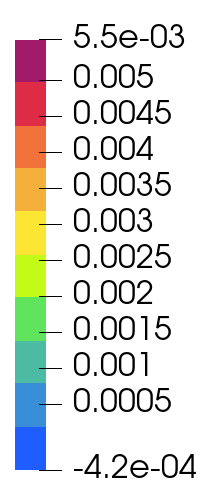} \\
    \textrm{(a)} & \textrm{(b)} & \\
    \multicolumn{1}{c}{\includegraphics[width=0.2\textwidth]{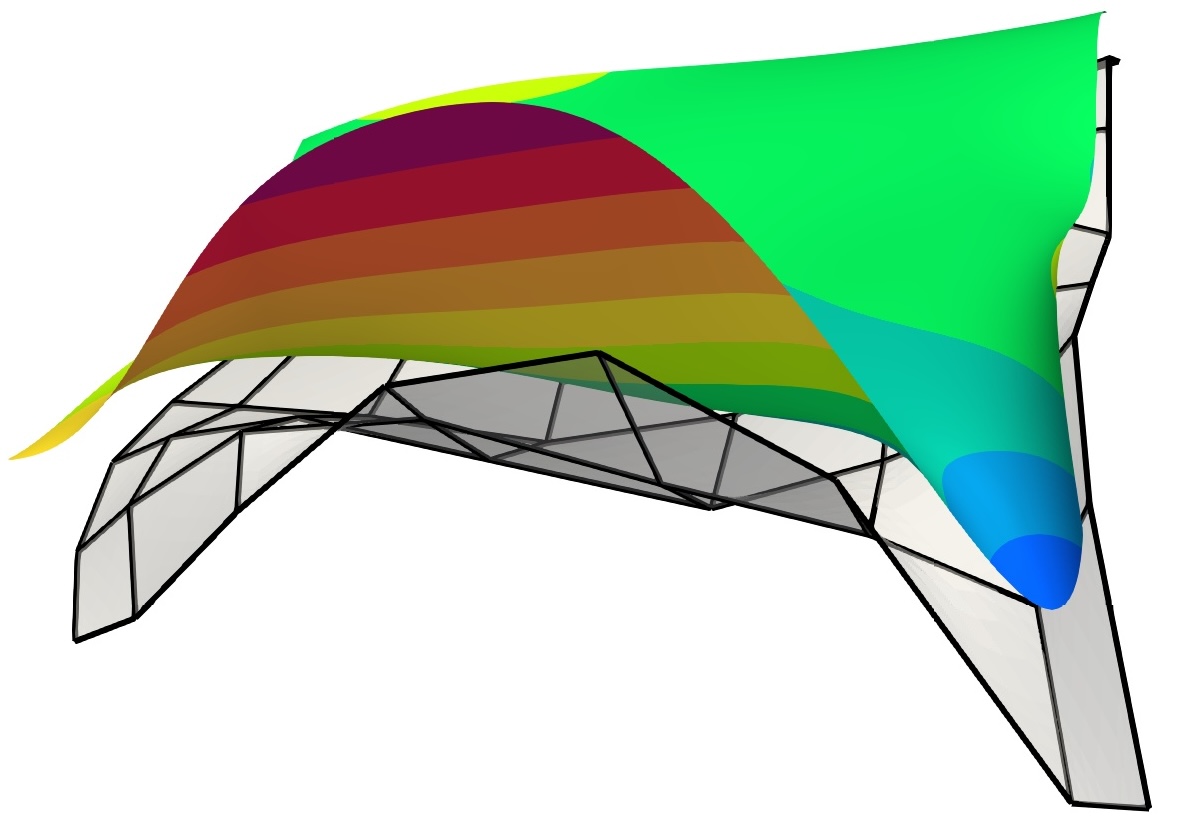}} &
    \multicolumn{2}{c}{\includegraphics[width=0.23\textwidth]{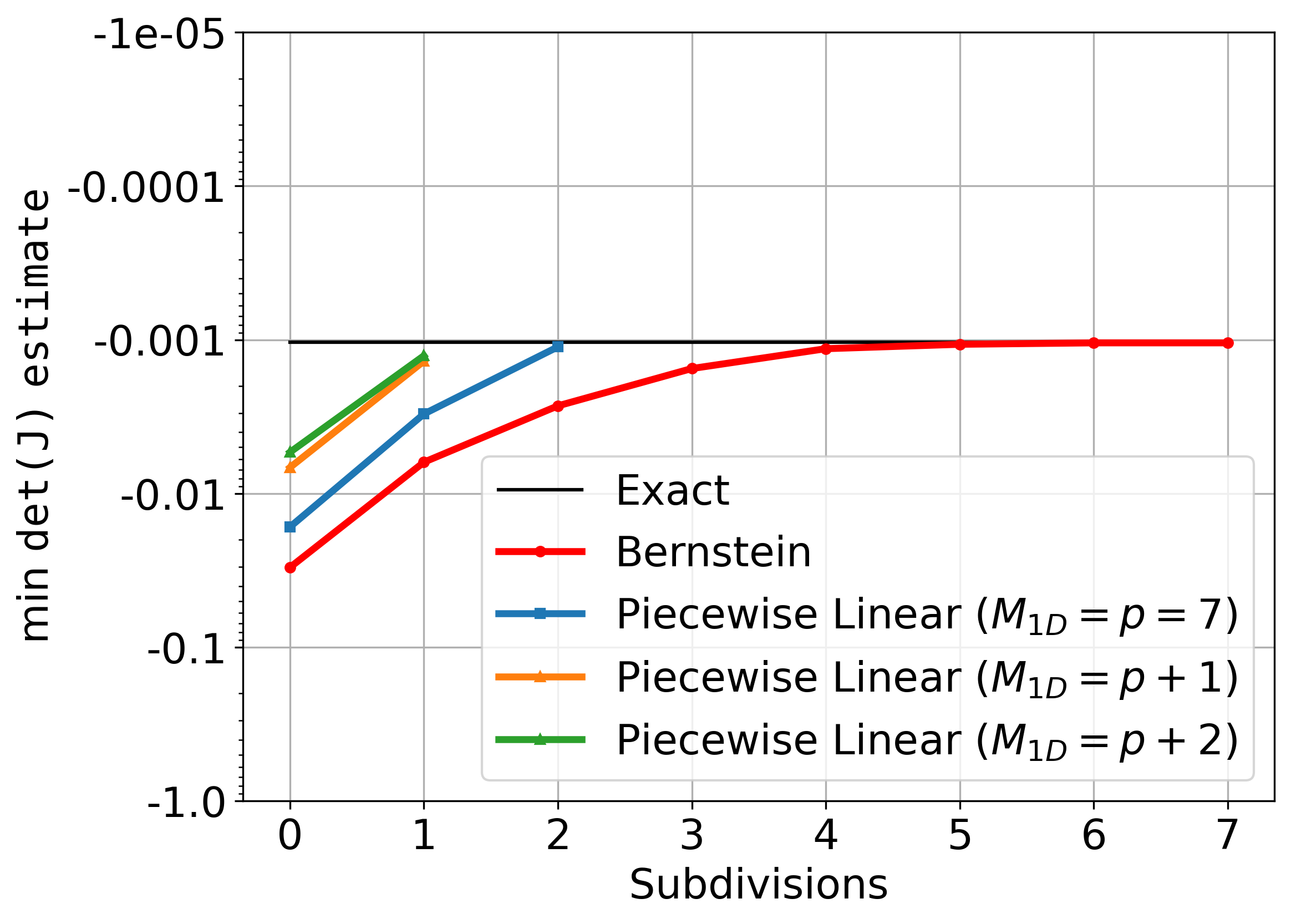}} \\
    \multicolumn{1}{c}{\textrm{(c)}} &
    \multicolumn{2}{c}{\textrm{(d)}} \\
\end{array}$
\end{center}
\vspace{-1mm}
\caption{(a) Inverted element from the blade example in Figure \ref{fig:blade_example}. (b) Jacobian determinant in the reference configuration. (c) Lower piecewise linear bound on the jacobian determinant with $8 \times 8$ control points. (d) Comparison of the lower bound using different number of control points and Bernstein polynomials with exact minimum.}
\label{fig:single_quad_detj}
\vspace{-3mm}
\end{figure}

Figure \ref{fig:single_quad_detj} shows the 4th-order inverted element from the blade example in Figure \ref{fig:blade_example} and the corresponding 7th-order Jacobian determinant on the reference element. The approach from \cite{mittal2025general,dzanic2025method} is used to obtain the lower bound on the Jacobian determinant on a control net of $8 \times 8$ points, as illustrated in Figure \ref{fig:single_quad_detj}(c). Finally, Figure \ref{fig:single_quad_detj}(d) illustrates how the lower bound converges towards the exact minimum using Bernstein polynomials and piecewise linear bounds for increasing number of control points. For the piecewise linear bounds, we recursively bound the portion of the element where the upper and lower bounds have opposite signs, until the upper bound becomes negative at any control point; this condition is sufficient to guarantee that the element is inverted.
With Bernstein bases, we do not have a lower and upper bound at each point, so we $h$-refine the element until the difference between the maximum and minimum coefficient of the bases is below a certain threshold.
\rev{As evident in Figure \ref{fig:single_quad_detj}(d), the piecewise linear bounds converge to the exact minimum significantly faster than Bernstein-based bounds, making the proposed approach more efficient at improving mesh quality and guaranteeing mesh validity.}

\subsubsection*{Mesh untangling}
When the initial mesh is invalid, we also use $\underline{\alpha}$ in shifted-barrier metrics:
\begin{eqnarray}
    \label{eq:shifted_barrier_bounds}
    \mu(T) = \frac{\tilde{\mu}(T)}{2\big(\tau - \tau_b\big)}, \\
    \tau_{b} = \begin{cases}
    \underline{\tau} \equiv \underline{\alpha}\,{\omega^{-1}} - \epsilon & \text{if } \underline{\alpha} \leq 0 \\
    0 & \text{otherwise}
    \end{cases}
\end{eqnarray}
where $\omega$ refers to the determinant of the target matrix $W$, which is constant throughout the domain for the problems considered in this manuscript.
The $\epsilon$ offset is used here because the bounding approach from \cite{mittal2025general} can capture the exact minimum if it is located at an element vertex.
Unlike in \eqref{eq:shifted_barrier_qp}, the barrier $\tau_b$ \eqref{eq:shifted_barrier_bounds} is  guaranteed to be below $\tau_{\tt min}$.
In spirit, this approach is similar to \cite{toulorge2013robust} where the bound on the Jacobian determinant, obtained via Bernstein bases, is used to set a barrier in the functional for $r$-adaptivity.

\subsubsection*{\rev{Computational Cost}}
\rev{
Assuming the nodal position function is defined
using $N_p$ bases, the computational cost to evaluate the Jacobian determinant at $M_q$ quadrature points is $\mathcal{O}(N_p \cdot M_q)$. We have typically used $M_q = 2 N_p$ in our prior work.
Similarly, the cost of representing the Jacobian determinant as
a high-order function with $N_{\alpha}$ nodes and Lagrange interpolants is $\mathcal{O}(N_p \cdot N_{\alpha})$.
These nodal function values can then be used to compute the Bernstein coefficients using a change-of-bases transformation matrix, which costs $\mathcal{O}(N_{\alpha}^{2})$. They can also be used to compute the piecewise linear bounding function on $M$ control points, which costs $\mathcal{O}(N_{\alpha} \cdot M)$, with $M \approx N_{\alpha}$.
Thus, the use of piecewise linear bounding functions
is more expensive than the \emph{non-rigorous} discrete sampling strategy but equivalent to the \emph{conservative} Bernstein-based approach. Note that in addition to generally returning looser bounds, the Bernstein-based approach also suffers from poor conditioning of the change-of-bases transformation matrix, especially for higher-order polynomials.}

\subsection{Tangential relaxation on curved surfaces} \label{subsec:tangentialrelaxation}

We use the high-order mesh itself as a proxy for geometry
to enable tangential relaxation without relying on an external CAD model.
Our approach allows all nodes, including those on curved boundaries, to move freely in response to the mesh optimization objective, and then restores the boundary nodes back to the original mesh surface using closest point projection. The resulting displacement is smoothly propagated into the interior via a Laplace solve.

\rev{
Since the initial mesh serves as the geometric surrogate, the accuracy of the final boundary representation is limited by the fidelity of the input mesh. This is a trade-off for the CAD-free capability and we assume that the initial
mesh captures the geometry of the domain with sufficient accuracy. }

\subsubsection*{\rev{Closest Point Projection}}
Closest point projection is enabled by leveraging the framework of \cite{mittal2025general} designed for fast and robust evaluation of high-order functions on curved meshes for any given point $\bx^*$ in physical space.
This entails computing the overlapping element ($\Omega^{e^*}$) and corresponding reference coordinates ($\bbx^*$) for $\bx^*$, such that the desired finite element function can be evaluated at that location (e.g., using \eqref{eq:u}).
The framework in \cite{mittal2025general} does this efficiently by first using a combination of local and global map to identify candidate elements and corresponding mpi ranks that could overlap $\bx^*$. Next, element-wise bounding boxes are used to further narrow down the list of candidate elements. Finally, a nonlinear minimization problem is solved in each of these elements to identify the element that contains the point and compute the corresponding reference-space coordinates:
\begin{equation}
    \boldsymbol{\bbx}^* = \arg\min_{\boldsymbol{\bbx} \in \bar{\Omega}^e}
    \left\| \bx^* - \bx(\bbx) \right\|^2,
    \label{eq:closest-point}
\end{equation}

\subsubsection*{\rev{Displacement Blending}}
During $r$-adaptivity, we minimize the TMOP objective \eqref{eq:F_full} while allowing the boundary nodes on curved boundaries ($\partial \Omega_T$) to move freely based on the solver update \eqref{eq:r_adaptivity_solve} and obtain tentative nodal positions, $\tilde{\bx}_{k+1}$.
Then, for each boundary node, we solve \eqref{eq:closest-point} with the surface mesh extracted from the initial mesh, which naturally results in closest point projection. The updated positions $\breve{\bx}_{k+1}$ contain projected positions for nodes on $\partial\Omega_{T}$ and unchanged positions for the interior nodes.
The boundary nodes cannot be independently moved to the surface so we blend the displacements,
\begin{equation}
    \Delta \breve{\bx}_{k+1} = \breve{\bx}_{k+1} - \tilde{\bx}_{k+1},
    \label{eq:boundary-displacement}
\end{equation}
into the interior of the domain by solving the Laplace problem:
\begin{equation}
\begin{aligned}
    -\nabla^2 (\Delta {\bx}_{k+1}) &= \mathbf{0} && \text{in } \Omega, \\
    \Delta {\bx}_{k+1} &= \Delta \breve{\bx}_{k+1} && \text{on } \partial\Omega_{T}, \\
    \Delta {\bx}_{k+1} &= \mathbf{0} && \text{on } \partial\Omega \setminus \partial\Omega_{T},
\end{aligned}
\label{eq:laplace-extension-full}
\end{equation}
where $\partial\Omega \setminus \partial\Omega_{T}$ denotes boundaries where nodal positions are fixed.
Finally, the updated positions are obtained by adding the blended displacement back to positions produced by the Newton update:
\begin{eqnarray}
    \label{eq:x_update_blend}
    \bx_{k+1} &=& \tilde{\bx}_{k+1} + \Delta {\bx}_{k+1}.
\end{eqnarray}

\begin{figure*}[tb]
\begin{center}
$\begin{array}{cccc}
    \includegraphics[width=0.2\textwidth]{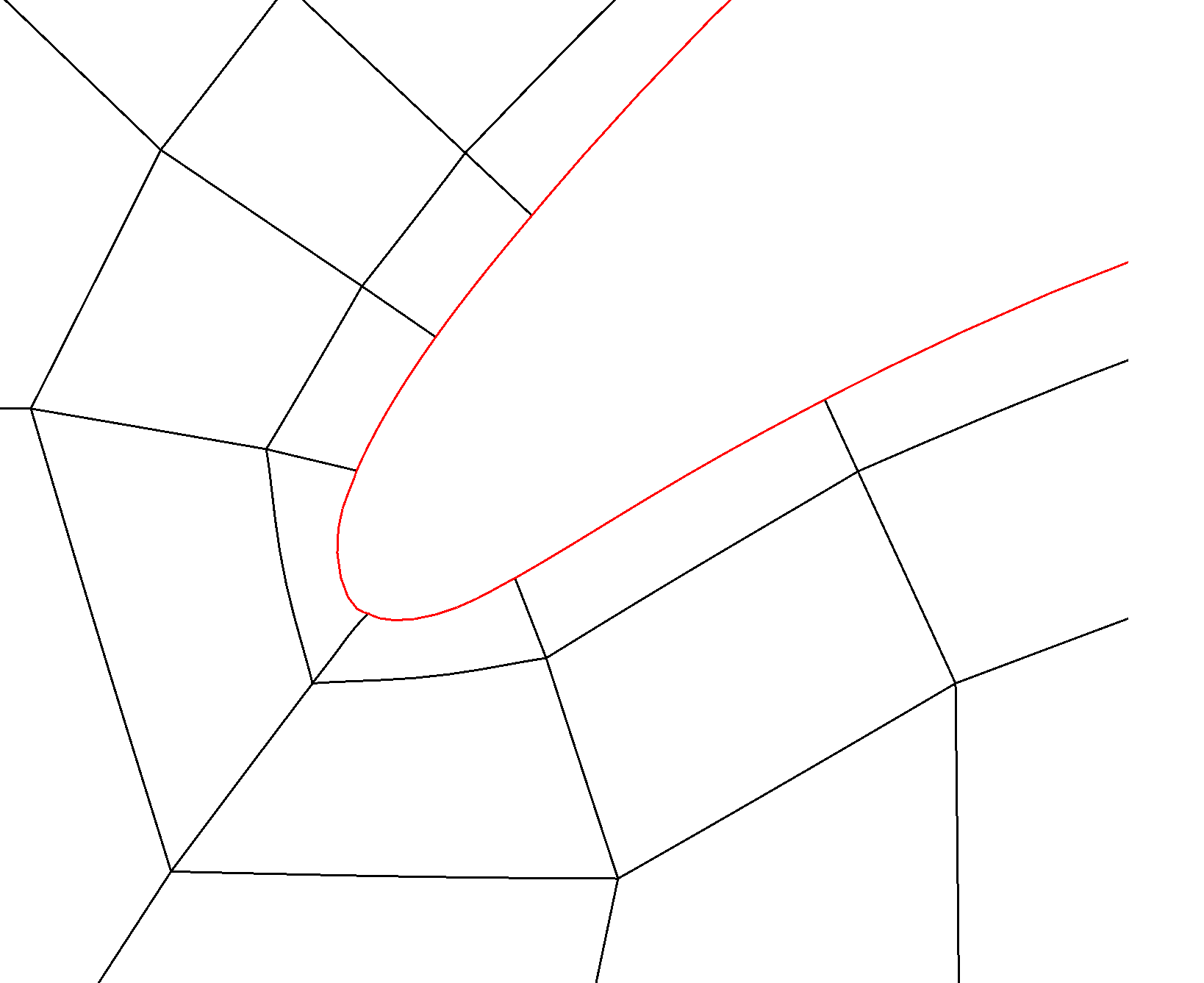} &
    \includegraphics[width=0.2\textwidth]{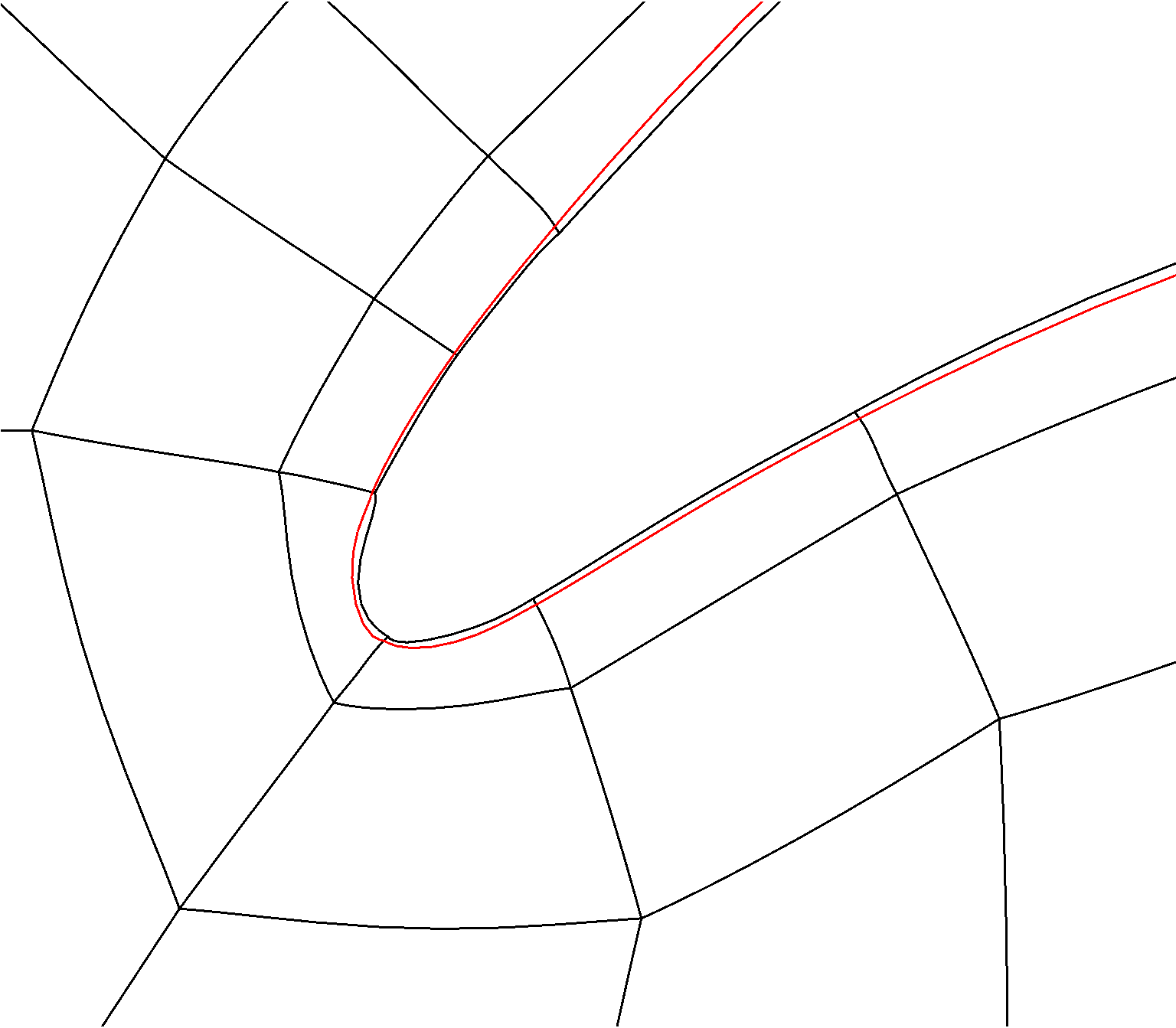} &
    \includegraphics[width=0.2\textwidth]{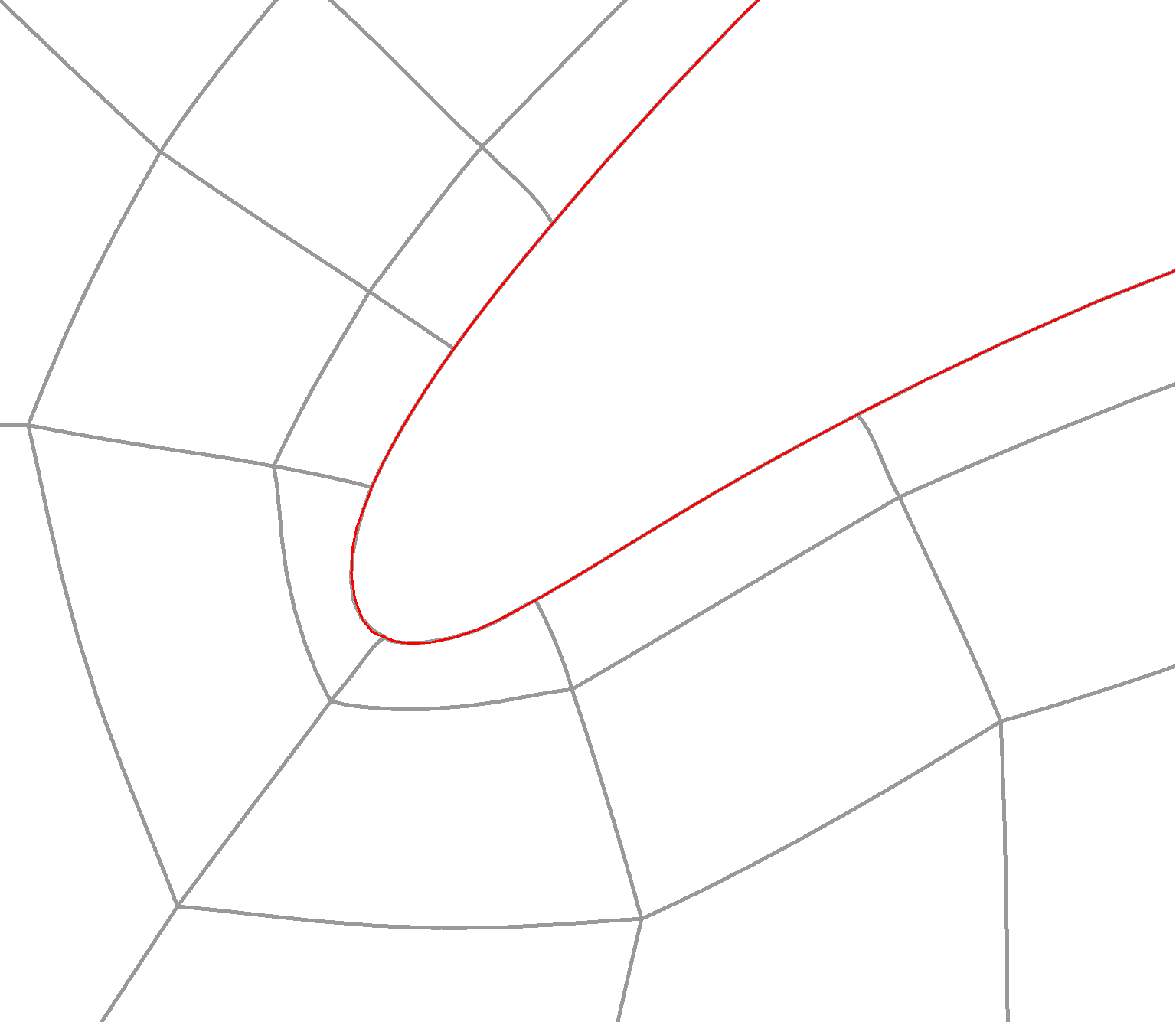} &
    \includegraphics[width=0.2\textwidth]{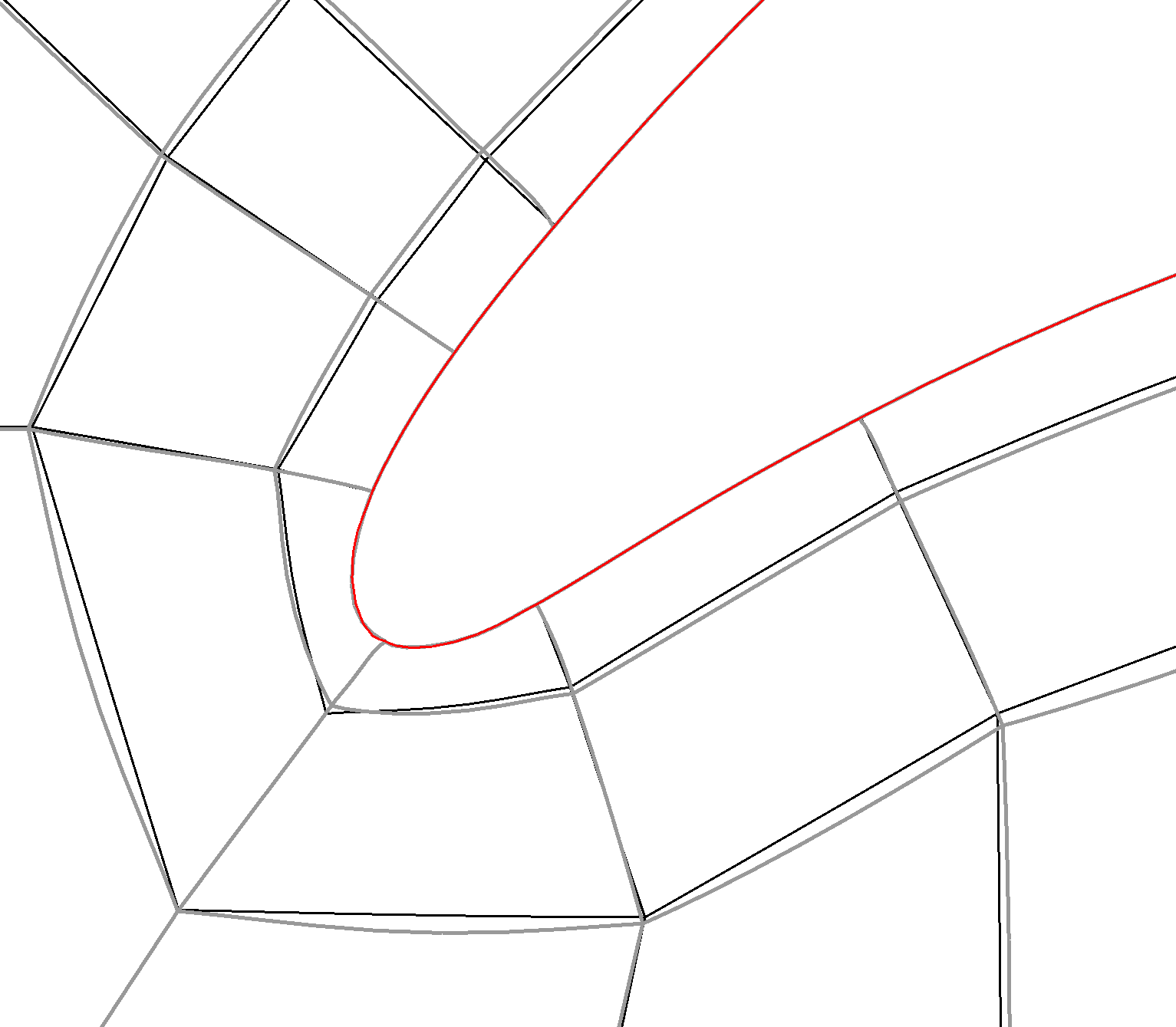} \\
    \textrm{(a)} & \textrm{(b)} & \textrm{(c)} & \textrm{(d)}
\end{array}$
\end{center}
\caption{Impact of tangential relaxation near the leading edge of the blade after one $r$-adaptivity iteration. The blade surface for tangential relaxation is shown in red.
(a) Initial mesh.
(b) Optimized mesh with boundary nodes moved due to the solver update.
(c) Boundary nodes projected back to the original mesh surface and corresponding displacement blended into the interior.
(d) Final mesh (gray) and input mesh (black) shown together for comparison.}
\label{fig:blade_tangential_relaxation}
\end{figure*}

This closest point projection algorithm is incorporated in the line-search step of our optimization framework. $\tilde{\bx}_{k+1}$ is the Newton update in \eqref{eq:r_adaptivity_solve}, and is used for closest point projection only if it forms a valid mesh, i.e. $\underline{\alpha}(\tilde{\bx}_{k+1}) > 0$. After the closest point projection and displacement blending, the updated positions $\bx_{k+1}$ are accepted only if they form a valid mesh ($\underline{\alpha}({\bx}_{k+1}) > 0$), and  satisfy the energy \eqref{eq:line_search_1} and norm of the gradient \eqref{eq:line_search_2} constraint.

We have observed that with the Newton's method, the mesh motion can be quite drastic with the boundary nodes unconstrained. This can be especially problematic near concave surfaces. Using a gradient-based solver such as the Broyden–Fletcher–Goldfarb–Shanno (BFGS) method helps mitigate this issue, and enables robust and accurate tangential relaxation while preserving the geometric integrity of the original boundary.

From an implementation perspective, the input mesh is preprocessed to associate each boundary node marked for tangential relaxation with a geometric entity: a corner, edge (in 2D and 3D), or surface (in 3D). The corresponding edge or surface mesh is also extracted.
This preprocessing step is critical for maintaining geometric accuracy. We perform the node-to-geometry association by inspecting the \emph{attributes} of incident edges and faces. \emph{Attributes} are unique identifiers commonly used in meshing and finite element workflows to define boundary conditions.
In 2D, a node with a single attribute is associated with the corresponding edge and is allowed to move tangentially. A node with two distinct attributes is inferred to lie at a geometric corner and is held fixed during $r$-adaptivity. This approach assumes that adjacent mesh boundaries are tagged with different attributes. Similarly in 3D, nodes with one, two, or three distinct attributes are respectively classified as lying on a surface, an edge, or a corner.

Figure \ref{fig:blade_tangential_relaxation} illustrates this tangential relaxation near the leading edge of the blade. The initial mesh is shown in Figure \ref{fig:blade_tangential_relaxation}(a). After 1 iteration of $r$-adaptivity, the mesh nodes (including the boundary) move away from the initial mesh surface, as shown in Figure \ref{fig:blade_tangential_relaxation}(b). The boundary nodes are projected back to the original mesh surface using \eqref{eq:closest-point} and the displacement is blended into the interior, as shown in Figure \ref{fig:blade_tangential_relaxation}(c). The final mesh and input mesh are shown together in Figure \ref{fig:blade_tangential_relaxation}(d) for comparison.

\subsubsection*{\rev{Computational Cost}}
\rev{
The proposed tangential relaxation consists of
closest point projection followed by displacement blending at each backtracking trial during the line-search.
The closest point projection is performed independently for each boundary node and scales linearly with the number of boundary nodes; this is negligible compared to the global mesh optimization update owing to the efficient field evaluation framework \cite{mittal2025general}. The displacement
blending, however, entails a global Laplace solve (conjugate-gradient with preconditioned algebraic multigrid) that scales linearly with the number of nodes in the mesh.
Preliminary experiments indicate that the increase in computational cost due to displacement blending step of tangential relaxation may not be negligible as it can double the time spent in $r$-adaptivity.
However, the cost trade-off is justified when we consider scenarios where mesh quality improvement is severely restricted by fixed boundary nodes or when $r$-adaptivity is done once a preprocessing step before the actual simulation so its cost is easily amortized. Section \ref{subsec:results_ale} shows a 2D example where tangential relaxation makes a significant difference in element quality near the curved boundary.}

\section{Results}  \label{sec:results}

\subsection{Turbine blade}

\begin{figure*}[bt]
\begin{center}
$\begin{array}{ccc}
    \includegraphics[height=0.3\textwidth]{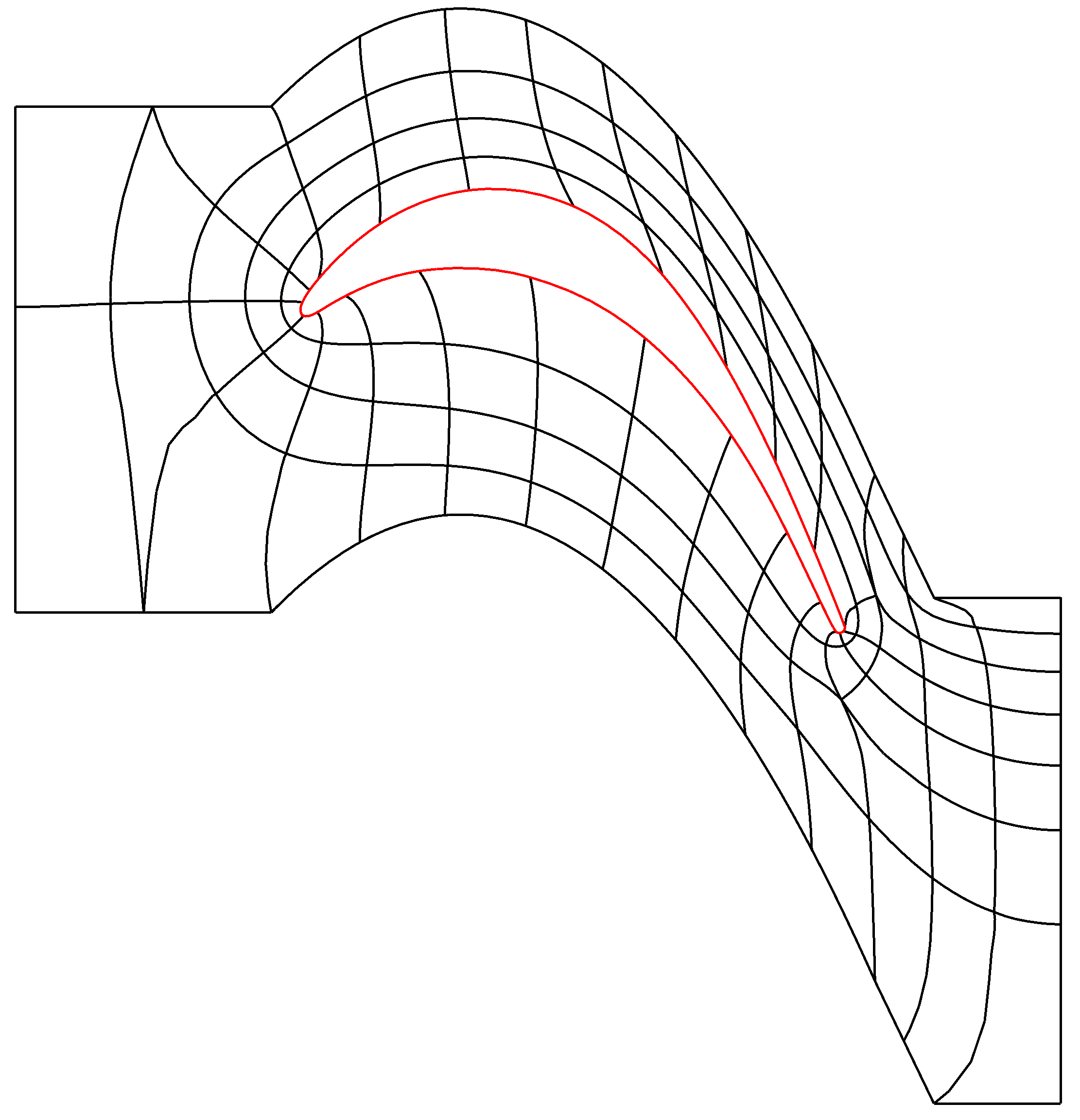} &
    \includegraphics[height=0.3\textwidth]{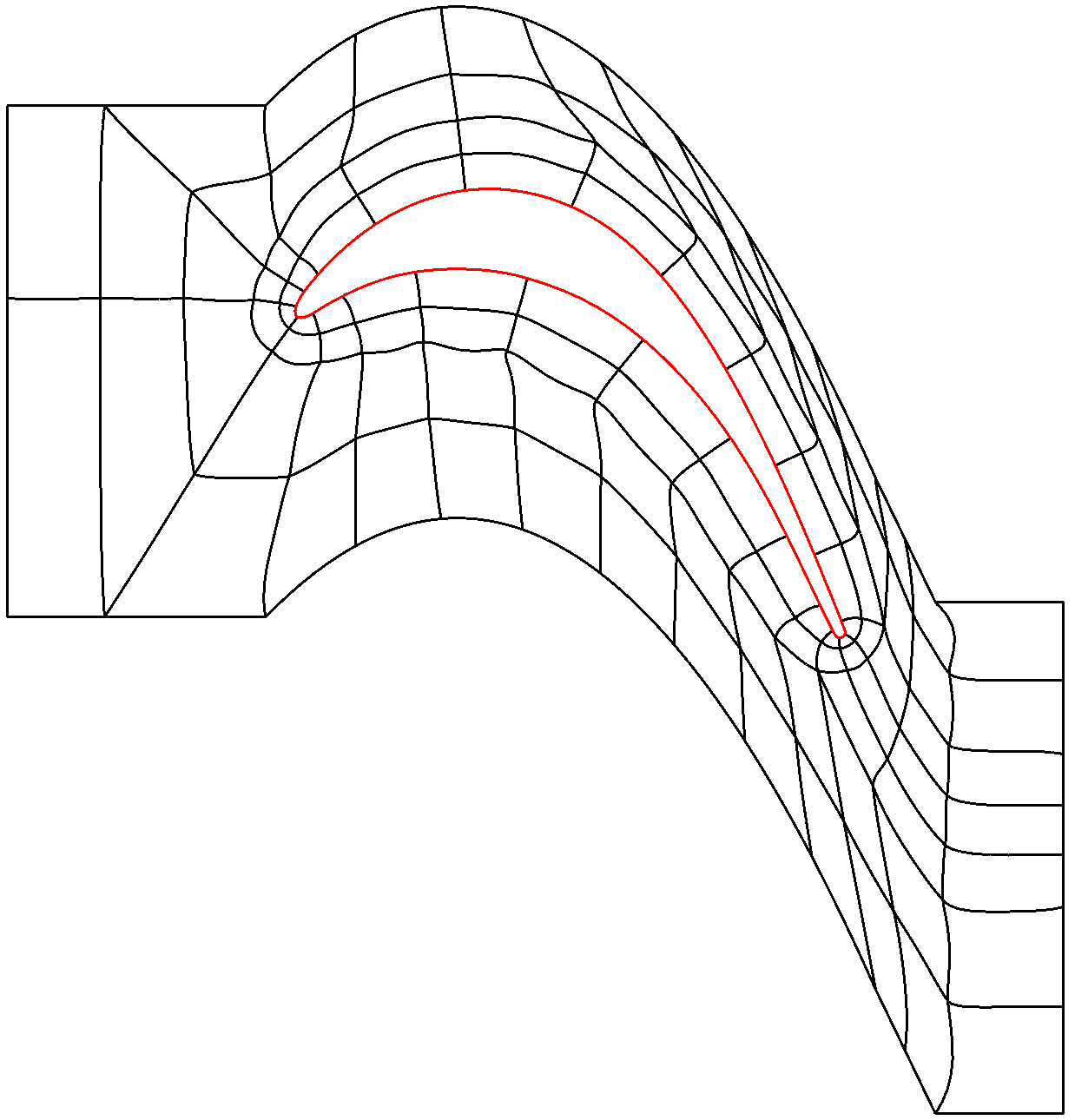} &
    \includegraphics[height=0.3\textwidth]{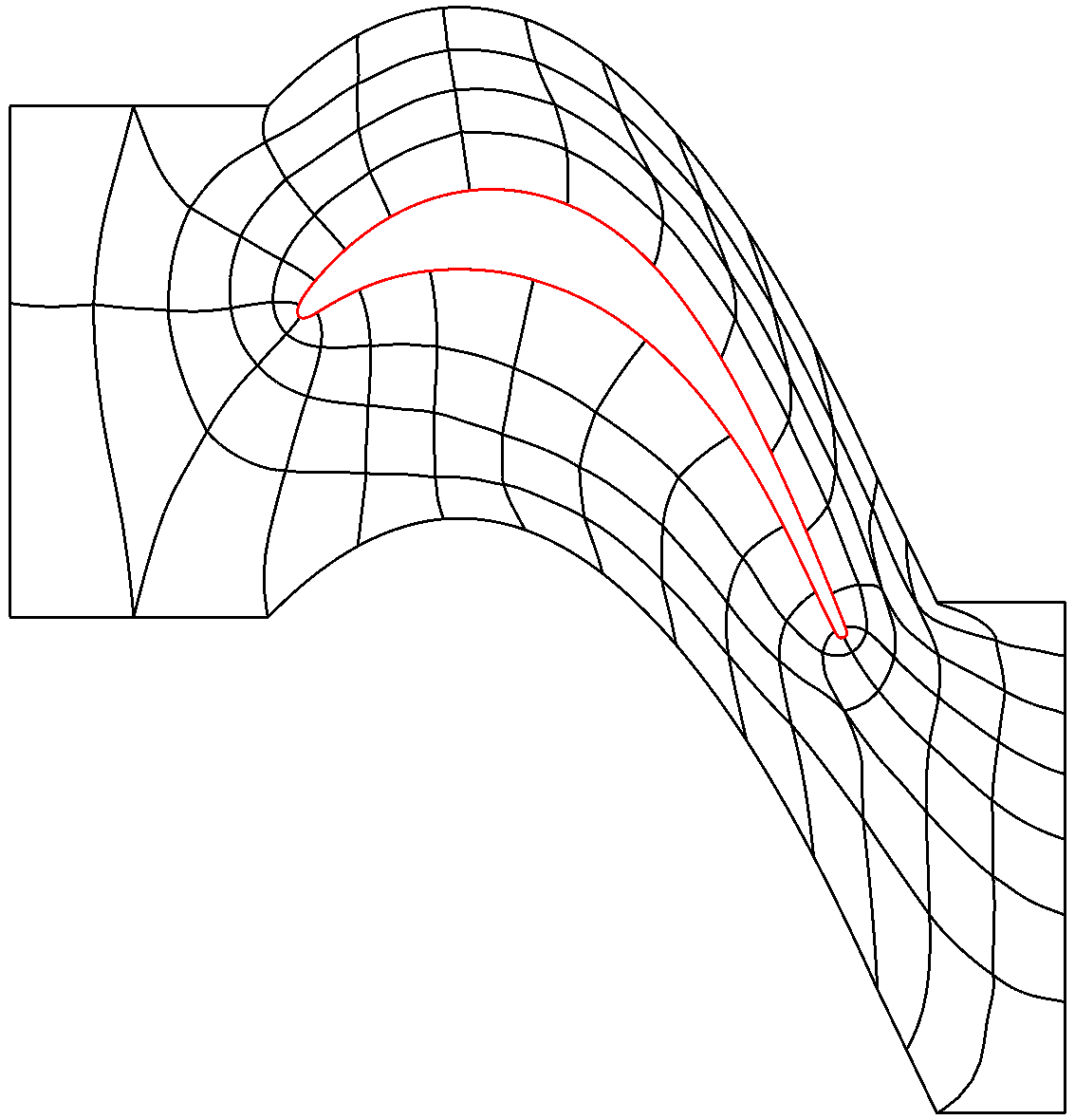} \\
    \textrm{(a)} & \textrm{(b)} & \textrm{(c)} \\
\end{array}$
\end{center}
\caption{Comparison of the optimized mesh for flow around a turbine blade. (a) Without tangential relaxation and mesh validity check. (b) Mesh optimization with tangential relaxation and guaranteed validity via bounds on Jacobian determinant stops prematurely as one of the element inverts
between quadrature points.
(c) Augmenting $r$-adaptivity with $q$-refinement enables significant improvement in mesh quality.}
\label{fig:blade_bounded}
\end{figure*}

We re-visit the 4th-order mesh for turbine blade example from Figure \ref{fig:blade_example}. The target matrix is still identity and a shape metric is used. A 10th-order accurate integration scheme with GLL quadrature points is employed for approximating the TMOP objective \eqref{eq:F_full}.
Figure \ref{fig:blade_bounded}(a) shows the optimized mesh from Figure \ref{fig:blade_example} where the mesh nodes are held fixed along the blade surface and the Jacobian determinant is only checked at quadrature points. Figure \ref{fig:blade_bounded}(b) shows the optimized mesh with the updates proposed in this manuscript. We observe that with tangential relaxation, the mesh elements are no longer skewed at the upper surface of the blade. However, the mesh optimization process stops prematurely due to the Jacobian determinant approaching 0 in the problematic element.
The TMOP objective is 170.8 for the initial mesh and reduces to 64.9 for the mesh in Figure \ref{fig:blade_bounded}(a). This mesh is not continuously valid however with $\ua(\bx)=-0.001$. For the mesh in Figure \ref{fig:blade_bounded}(b), the TMOP objective is 98.8 but the mesh is valid everywhere with $\ua(\bx)=10^{-10}$. In both cases, the minimum sampled Jacobian determinant is $10^{-4}$.

Increasing the integration order uniformly can help address the problem but will make the computational cost of $r$-adaptivity intractable. We thus employ $q$-refinement and use the ratio of the minimum sampled Jacobian determinant and its lower bound serves as an \emph{error estimator}; the quadrature order is increased only for elements where this ratio exceeds a user-defined threshold, i.e.,
\begin{eqnarray}
\label{eq:q_refinement}
\alpha_{\tt qp,min}(\bx)/\ua(\bx) > \epsilon_q.
\end{eqnarray}
The threshold $\epsilon_q$ is set to 5.0 by default.
This $q$-refinement is done after each $r$-adaptivity iteration, and the quadrature order is increased by the initial quadrature order.
Figure \ref{fig:blade_bounded}(c) shows the optimized mesh with $q$-refinement. The mesh is still valid everywhere ($\ua(\bx)=10^{-10}$) and
the final mesh is noticeably better with TMOP objective reduced to 74.2.

The example in Figure \ref{fig:blade_bounded}(c) uses an initial integration order of 10 and allows a maximum quadrature order of 400. Different maximum quadrature orders and refinement threshold result in different local minima, but consistently produces better meshes than without $q$-refinement. In future work, we will explore schemes where new quadrature points are inserted adaptively based on the local Jacobian determinant bounds rather than using a standard quadrature rule (e.g., GLL). We will also explore gradient-based solvers that are not as sensitive to the $q$-refinement parameters.

\subsection{Mesh untangling}

\begin{figure*}[bth]
\begin{center}
$\begin{array}{cccc}
    \includegraphics[width=0.22\textwidth]{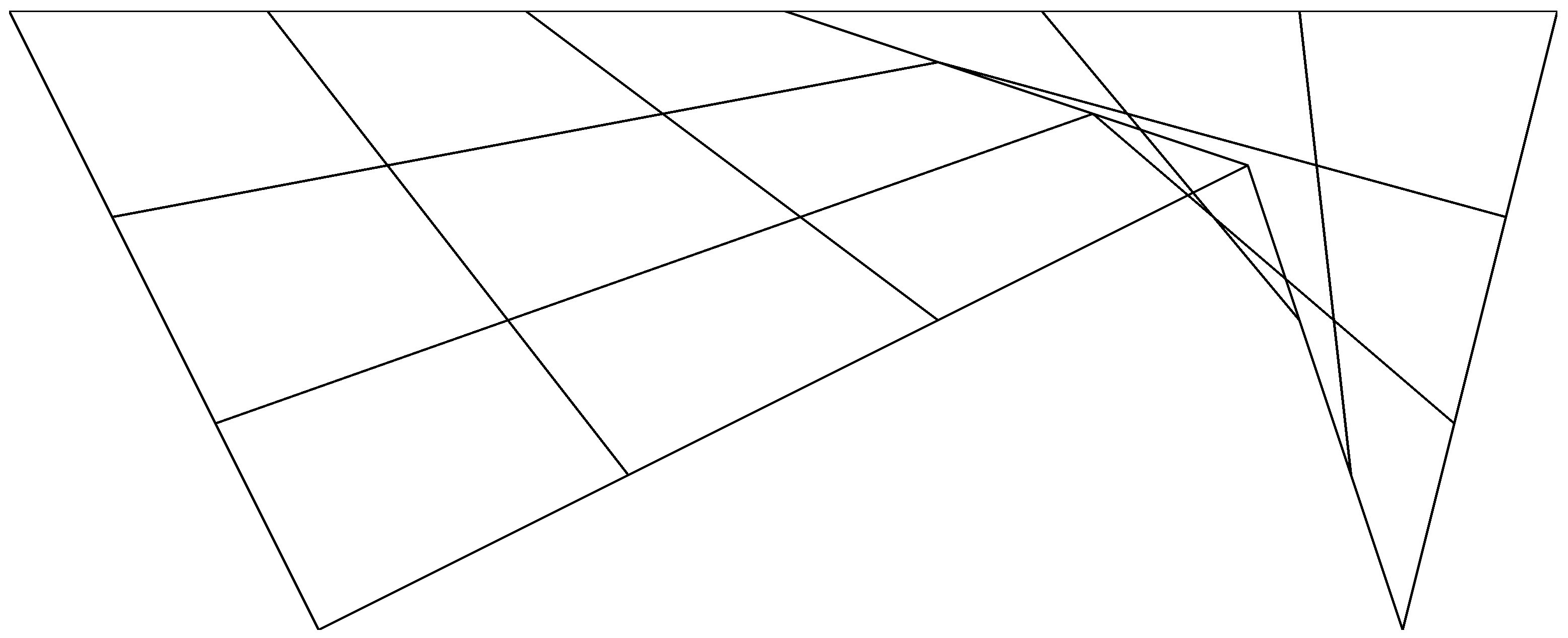} &
    \includegraphics[width=0.22\textwidth]{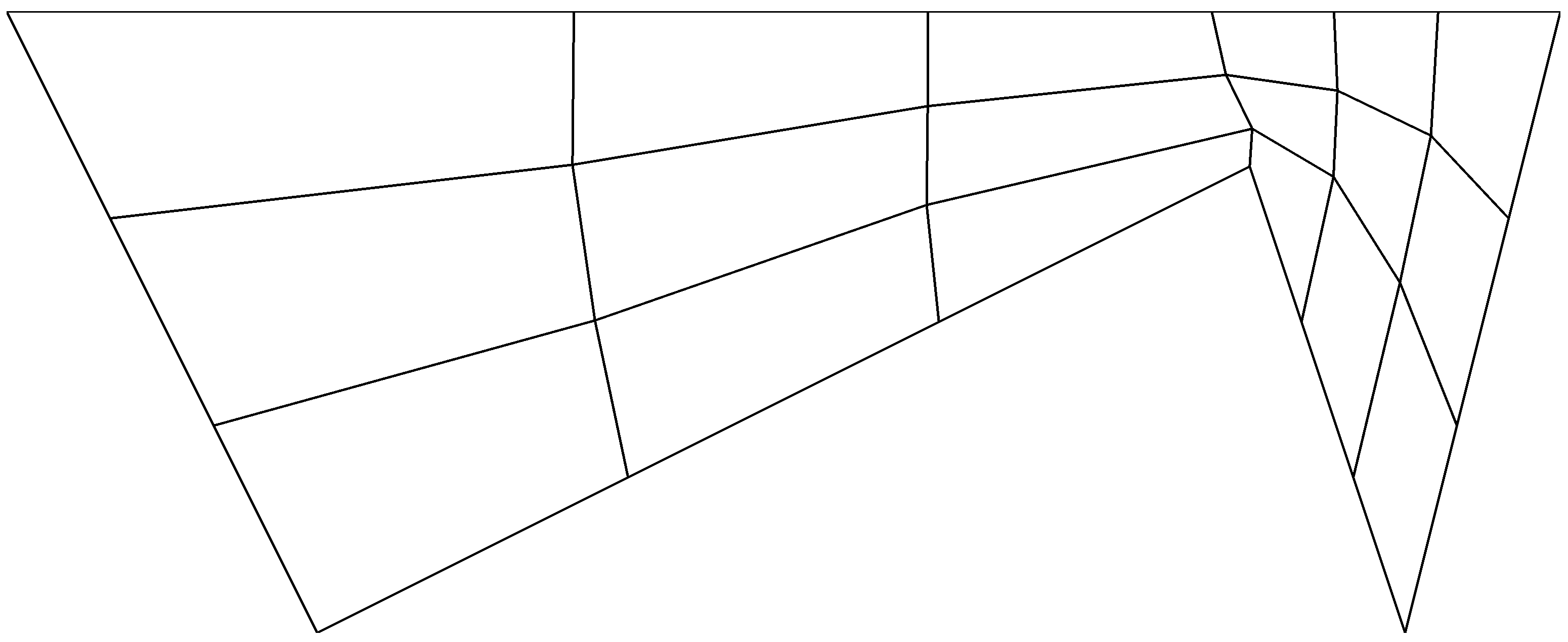} &
    \includegraphics[width=0.22\textwidth]{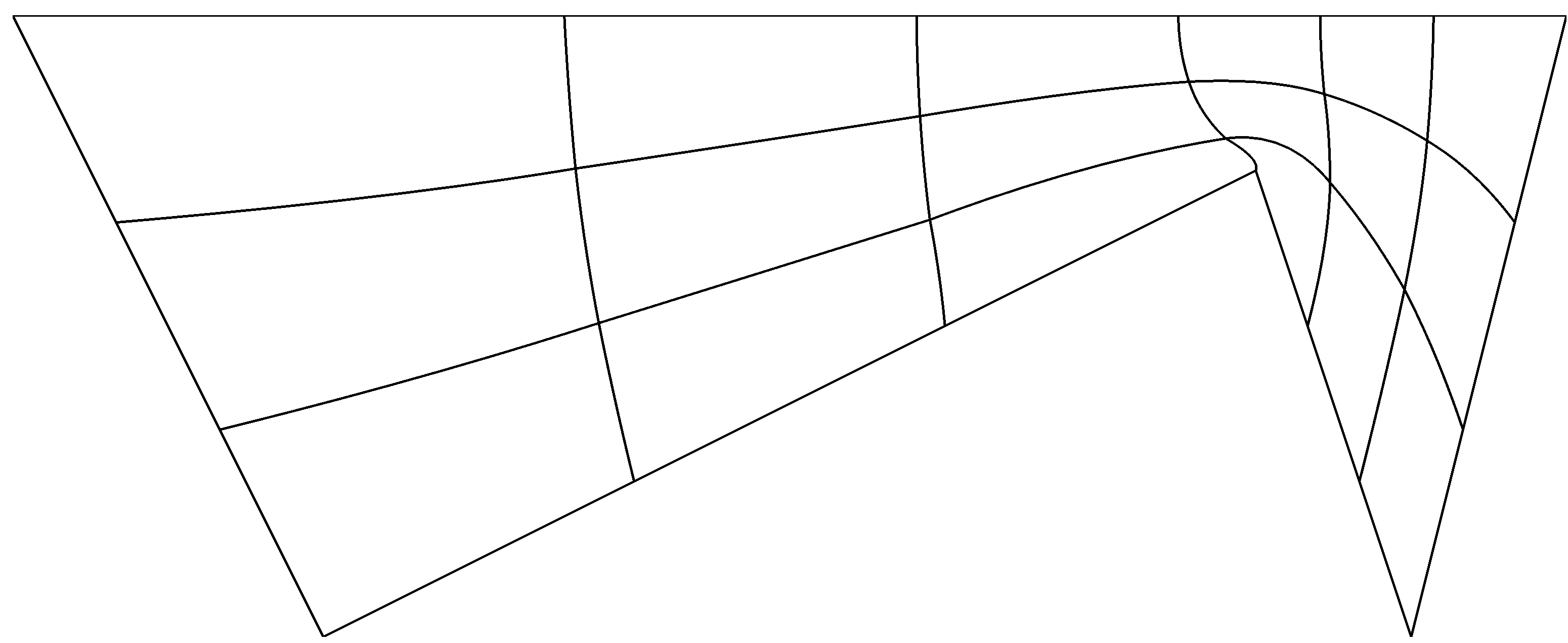} &
    \includegraphics[width=0.22\textwidth]{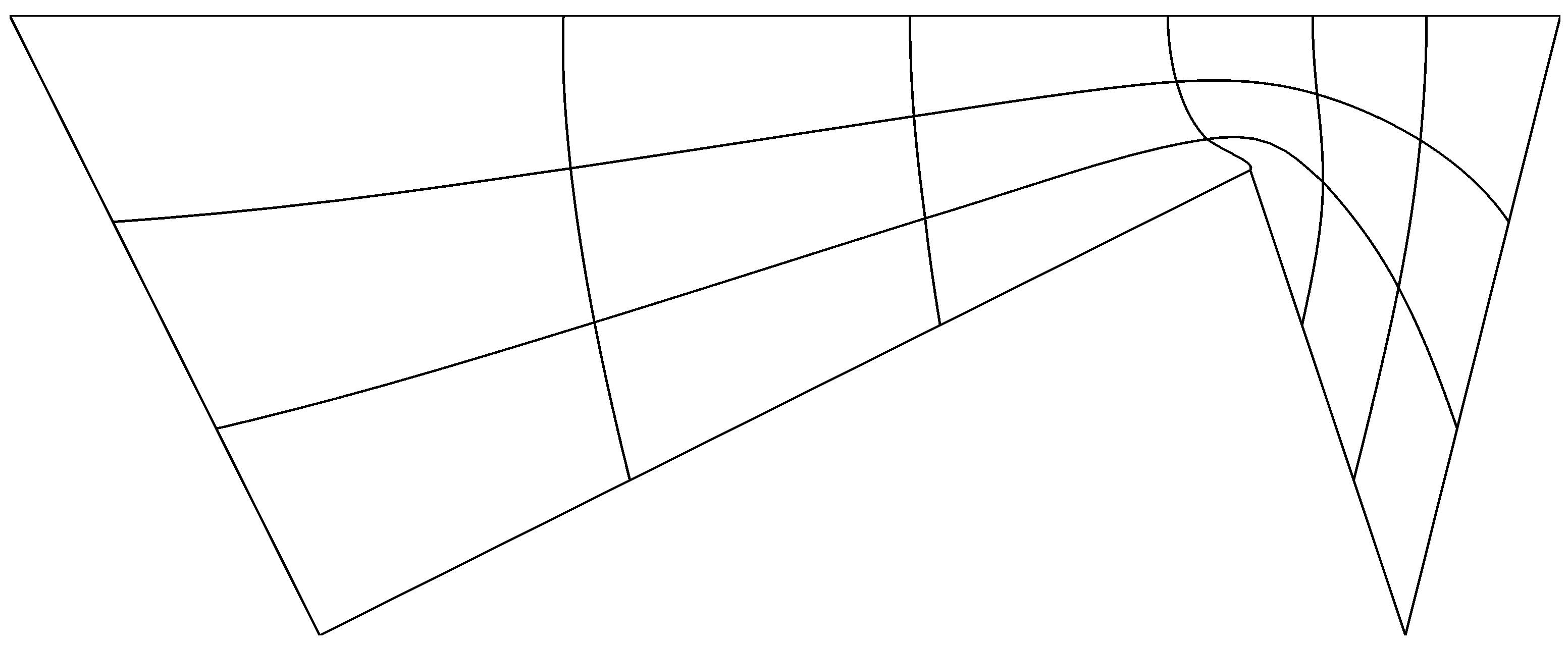} \\
    \textrm{(a) Initial mesh} & \textrm{(b)}\,\,p=1 &
    \textrm{(c)}\,\,p=2 & \textrm{(d)}\,\,p=3 \\
    \multicolumn{4}{c}{\textrm{(i) Mesh untangled using a shifted-barrier metric.}} \vspace{2mm}\\
    &
    \includegraphics[width=0.22\textwidth]{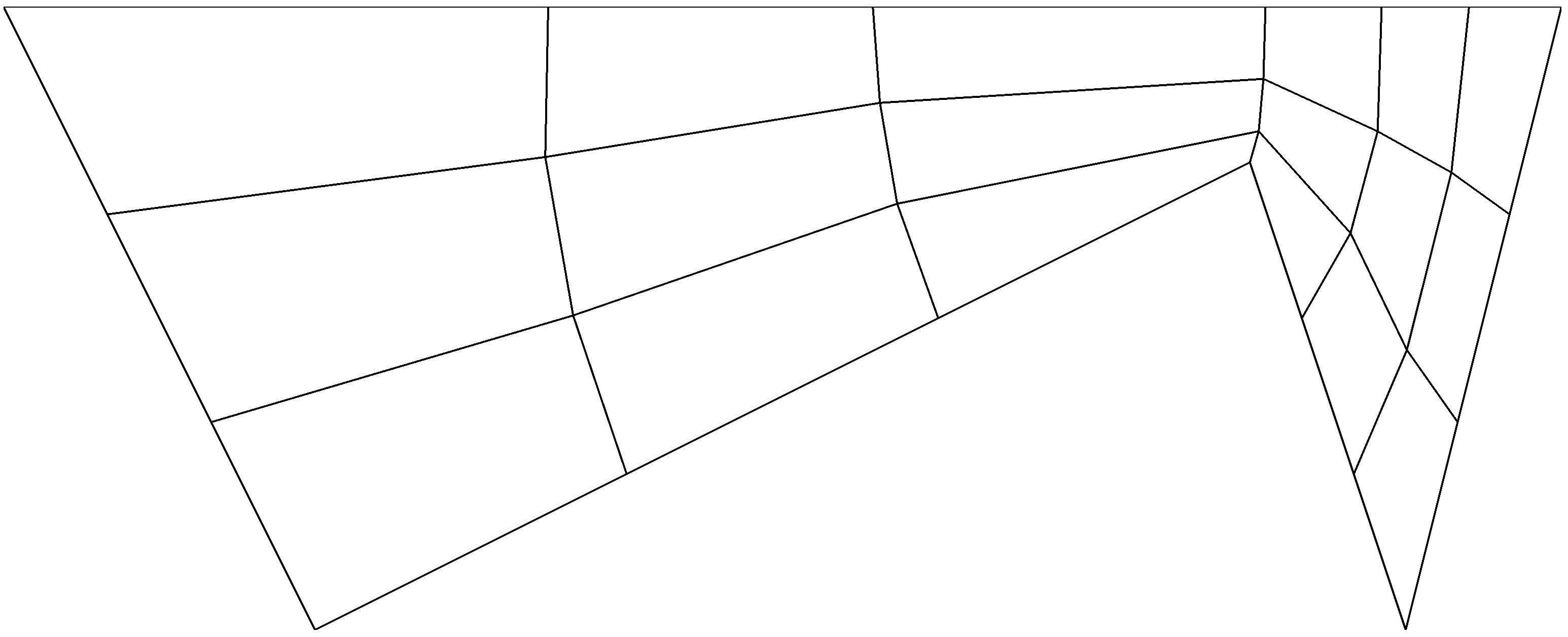} &
    \includegraphics[width=0.22\textwidth]{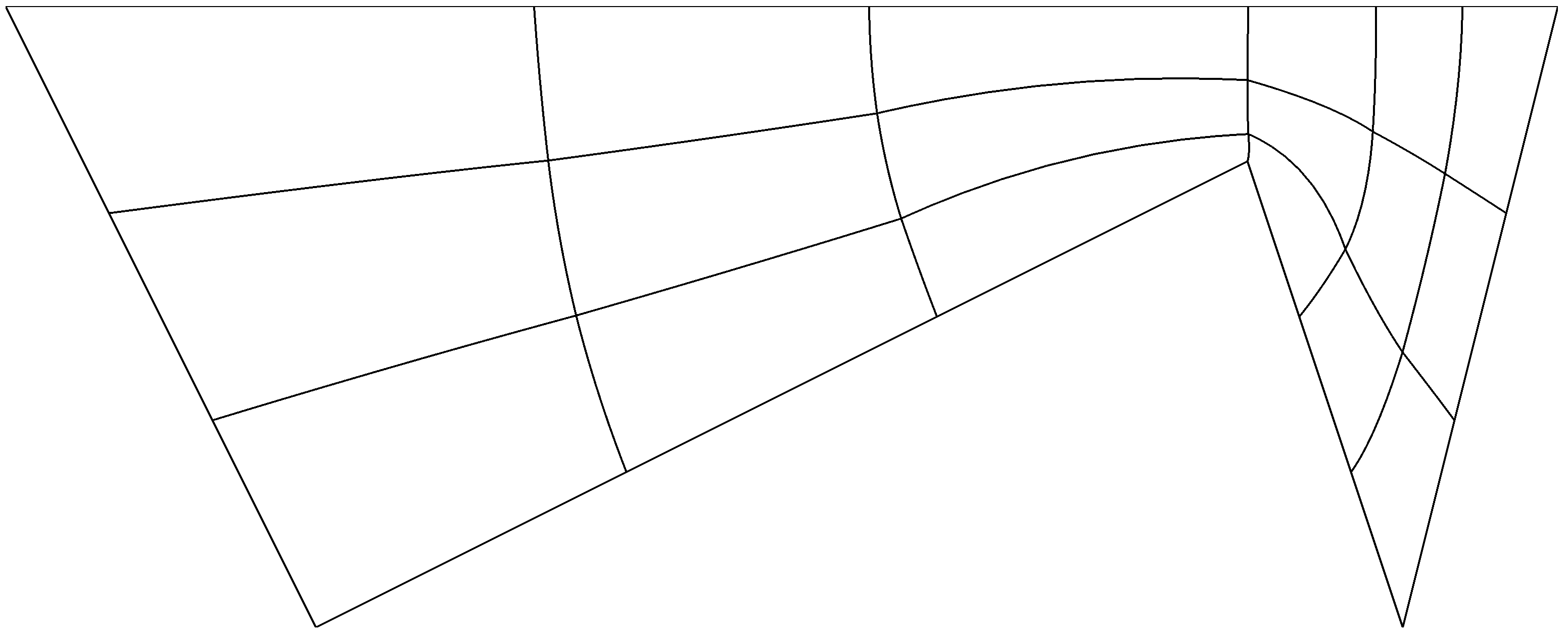} &
    \includegraphics[width=0.22\textwidth]{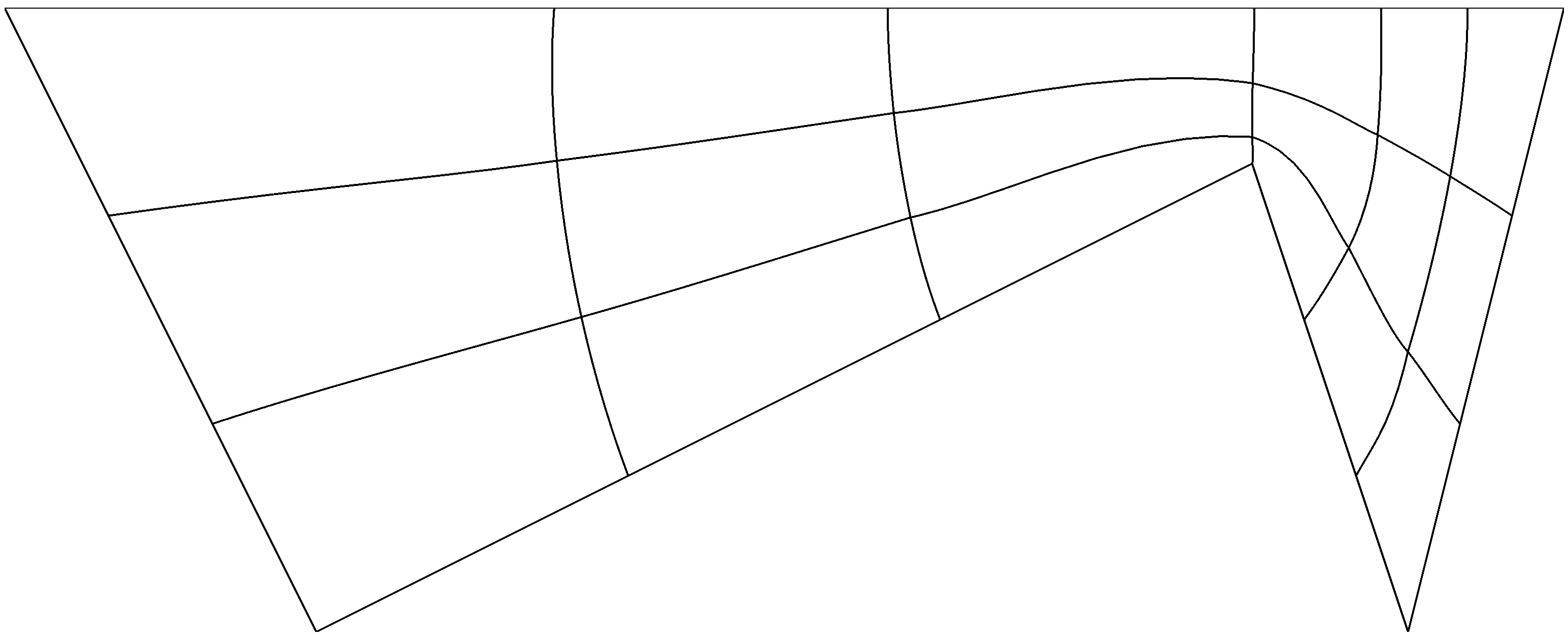} \\
    & \textrm{(e)}\,\,p=1 & \textrm{(f)}\,\,p=2 & \textrm{(g)}\,\,p=3 \\
    \multicolumn{4}{c}{\textrm{(ii) Untangled mesh improved for skewness.}}
\end{array}$
\end{center}
\caption{(i) Mesh untangling using a shifted-barrier metric with the lower bound on the Jacobian determinant as the barrier. (a) Initial mesh with inverted elements. Optimized mesh with polynomial order (b)  $p=1$, (c) $p=2$, and (d) $p=3$. The optimized mesh has positive lower bound on the Jacobian determinant in each case.
(ii) Further improvement of the untangled mesh for skewness using a convex combination of shape and skew metrics. (e) $p=1$, (f) $p=2$, and (g) $p=3$.
}
\label{fig:jagged}
\end{figure*}

We re-visit the mesh untangling example from \cite{knupp2022worst}, shown in Figure \ref{fig:jagged}(a). The mesh folds over itself, which results in inverting elements at the concave vertex. This poses a challenge for optimization-based smoothing, particularly at higher orders.

We optimize the mesh using a shifted-barrier approach with $\tilde{\mu}_4(T) = ||T||_F^2 - 2 \tau$ in \eqref{eq:shifted_barrier_bounds} with a target aspect ratio of unity and skewness of $\pi/2$. We use $\epsilon = 10^{-3}$ as the barrier offset.
Since tangential relaxation proposed in Section \ref{subsec:tangentialrelaxation} entails a PDE solve, and thus requires a valid mesh, it is disabled during mesh untangling unless the surface is aligned with the Cartesian axes.

\begin{figure}[tbh]
\begin{center}
$\begin{array}{c}
    \includegraphics[width=0.43\textwidth]{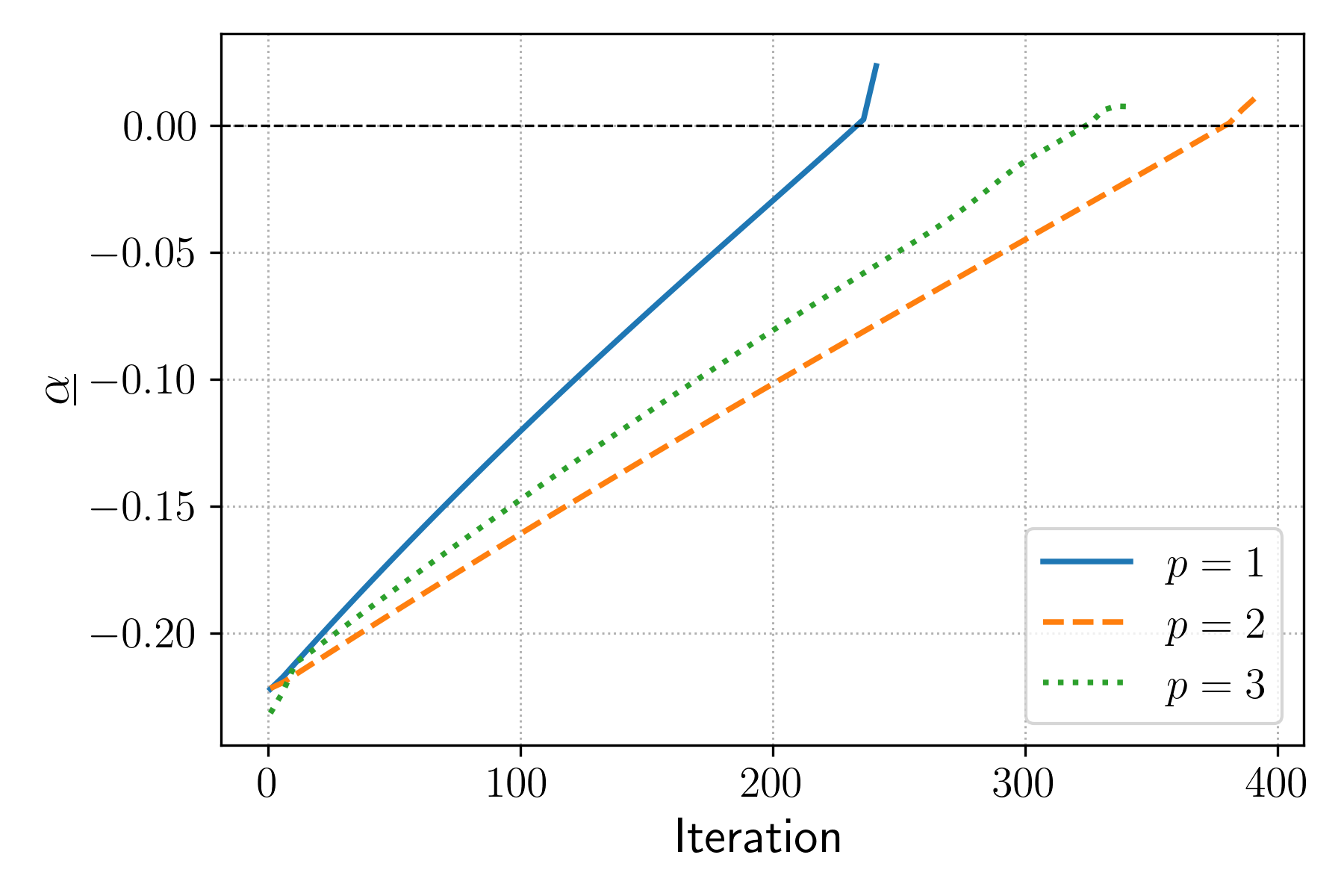}
\end{array}$
\end{center}
\caption{Evolution of the lower bound on Jacobian determinant during $r$-adaptivity with different polynomial orders.\vspace{-5mm}}
\label{fig:jagged_jac}
\end{figure}

Figure \ref{fig:jagged}(i) shows the optimized linear, quadratic, and cubic mesh. In each case, the mesh is successfully untangled. The lower bound on the Jacobian determinant of the initial mesh is -0.22 and reaches a positive value during $r$-adaptivity, as shown in Figure \ref{fig:jagged_jac}. Note that additional degrees of freedom in high-order elements allow the mesh to \emph{wrap} around the concave corner. The elements at this corner have poor skewness at the expense of aspect ratio trending towards unity. Since we cannot explicitly control the balance between aspect ratio and skewness in $\mu_{4,s}$, we further improve the untangled mesh using a convex combination of \emph{shape} metric and a purely \emph{skew} metric:
$\nu_{49,sq} = \gamma {\mu_{2,s}} + (1-\gamma)\nu_{50,q}$. Here, $\nu_{50,q}(A,W)=[1 - \cos(\phi_A - \phi_W)]/(\sin \phi_A \cdot \sin \phi_W)$ is the skew metric that depends on the local skewness angle $\phi$ of the current and target Jacobian of the transformation. The optimized meshes obtained with $\mu_{50}$ and $\gamma=0.4$ are shown in Figure \ref{fig:jagged}(ii), and as evident they have much better skewness than the elements of the untangled mesh in Figure \ref{fig:jagged}(i).

\subsection{Mesh relaxation for ALE Simulations}
\label{subsec:results_ale}

Arbitrary Lagrangian Eulerian (ALE) methods are widely used in multimaterial hydrodynamics to combine the advantages of Lagrangian and Eulerian frameworks. In the Lagrangian phase, the computational mesh moves with the material, capturing interfaces and conserving mass with high fidelity. However, strong deformations such as those arising in shock-driven flows can severely degrade mesh quality, necessitating periodic remeshing to maintain numerical accuracy and stability.

A typical ALE cycle consists of three stages: (i) Lagrangian solve on a moving high-order mesh, (ii) mesh optimization or rezoning to improve mesh quality while preserving geometric and material features, and (iii) remap of physical quantities to the new mesh configuration. The mesh optimization step is critical for ensuring that subsequent remap and hydrodynamic solves remain well-posed, especially in the presence of large deformations or evolving material boundaries. In this context, high-order mesh optimization plays a central role in enabling accurate, robust, and efficient ALE simulations.

In this section, we demonstrate the effectiveness of the proposed approach in improving mesh quality in the context of a 2D and 3D ALE simulation. In each example, we start with a high-order mesh that has been deformed by the Lagrangian solve, and then apply our $r$-adaptivity framework to optimize the mesh quality while preserving the geometric integrity of the domain boundaries.

\begin{figure*}[bth]
\begin{center}
$\begin{array}{cccc}
    \includegraphics[width=0.21\textwidth]{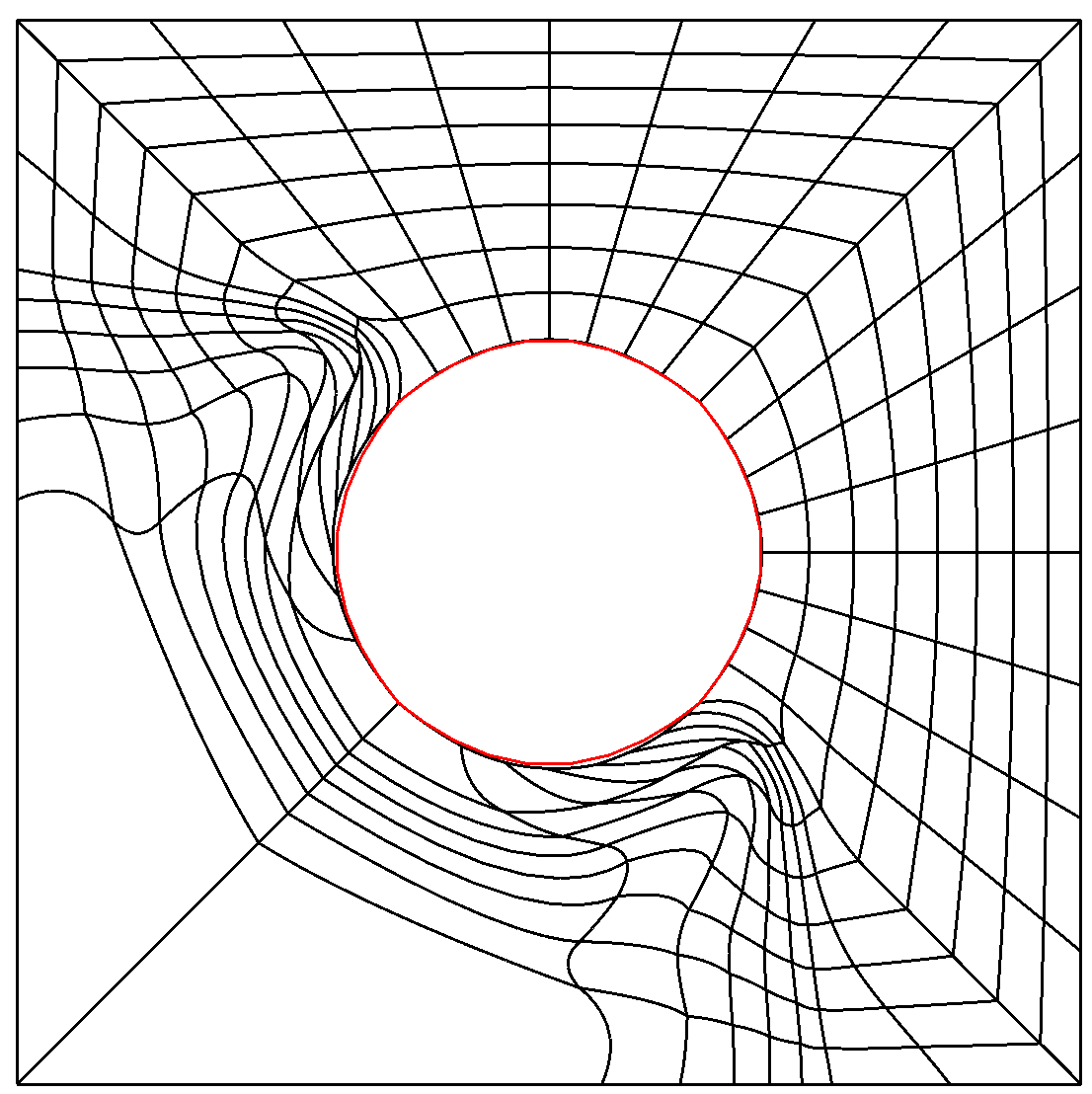} &
    \includegraphics[width=0.21\textwidth]{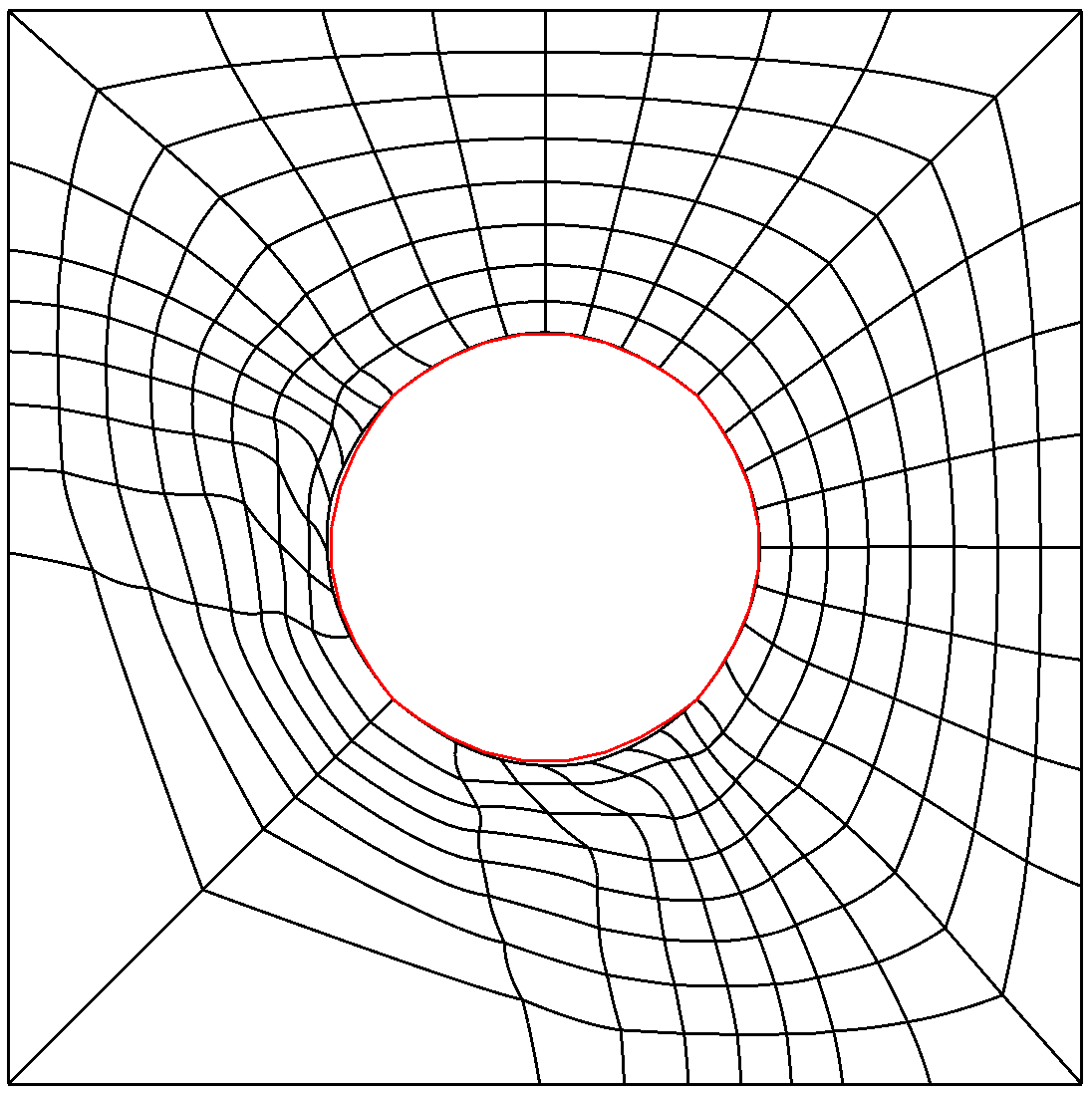} &
    \includegraphics[width=0.21\textwidth]{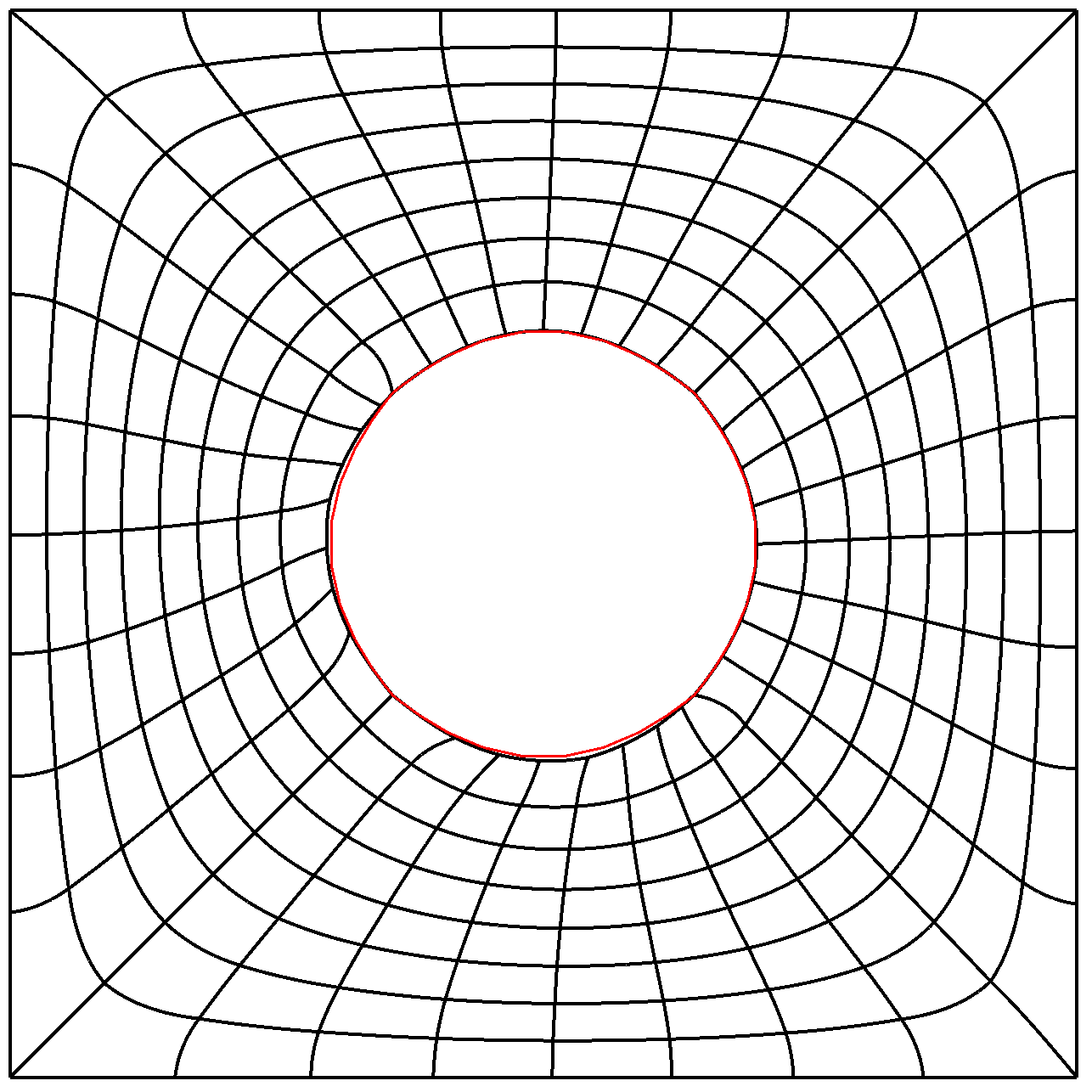} &
    \includegraphics[width=0.21\textwidth]{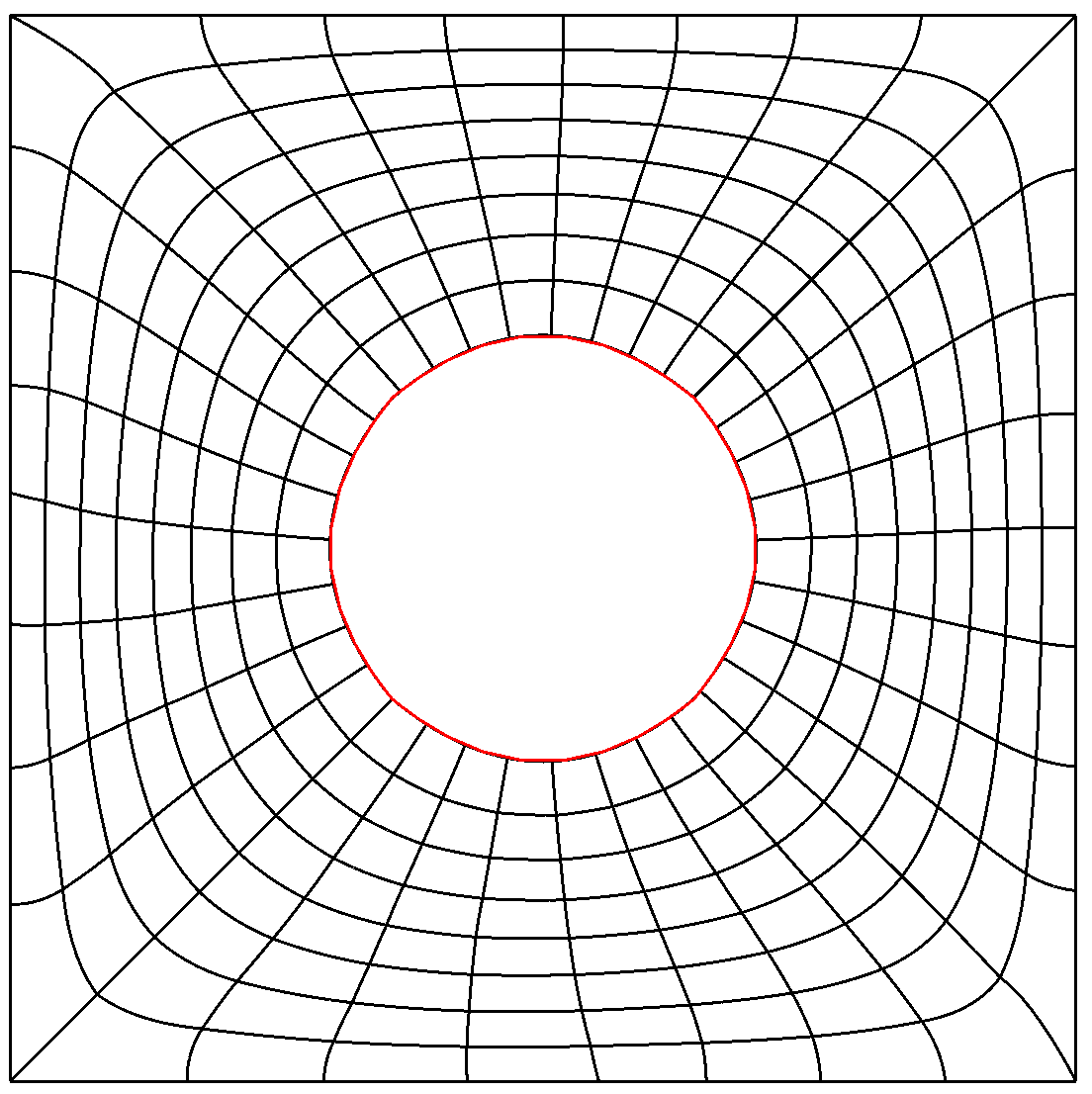} \\
    \textrm{(a) Initial mesh} &
    \textrm{(b) Untangled mesh} &
    \textrm{(c) Optimized mesh} &
    \textrm{(d) Optimized mesh} \\
    & & \rev{\textrm{(without tangential relaxation)}} &
    \textrm{(with tangential relaxation)}
\end{array}$
\end{center}
\caption{Mesh untangling using a shifted-barrier metric followed by $r$-adaptivity without and with tangential relaxation in the context of a 2D ALE simulation.}
\label{fig:2D_ale}
\end{figure*}

The 2D example is shown in Figure \ref{fig:2D_ale}.
The initial quadratic mesh is generated through a Lagrangian simulation of a
shock wave originating from the lower-left corner of the domain.
The wave collides with the circular obstacle and propagates around it;
see \cite{Dobrev2012, Nabil2024} for further details.
The deformed mesh is shown in Figure \ref{fig:2D_ale}(a). This mesh is valid at the quadrature points used during the Lagrangian phase, but has a negative Jacobian determinant between the quadrature points with a lower bound of $\ua(\bx) = -2 \times 10^{-4}$. The mesh is thus first untangled using a shifted-barrier metric ($\tilde{\mu}_4(T)$). This untangled mesh with $\ua(\bx) = 4 \times 10^{-4}$ is shown in Figure \ref{fig:2D_ale}(b). The boundary nodes on the curved boundaries are held fixed during mesh untangling.

Finally, the  mesh is optimized using an ideal shape \emph{target} with the shape metric $\mu_2(T)$. The optimized meshes without and with tangential
relaxation are shown in Figure \ref{fig:2D_ale}(c) and (d), respectively.
\rev{
The overall TMOP objective $\int \mu_2(T)$ is similar for both meshes as it decreases from 0.89 in the untangled mesh to 0.210 and 0.208 in the optimized meshes, respectively. The effectiveness of the tangential relaxation capability is evident when we consider the elements at the curved boundaries.
With the boundary nodes held fixed, the target circular boundary is not
recovered, $\ua(\bx) = 8 \times 10^{-4}$, and the skewness angle varies
from $\theta = 44^\circ$ to $\theta = 136^\circ$.
With tangential relaxation, the boundary nodes conform to the target
circular boundary, $\ua(\bx) = 2 \times 10^{-3}$, and the element edges
are closer to being orthogonal to the boundary with skewness angle
varying between $\theta = 82^\circ$ and $\theta = 97^\circ$.
Thus, the tangential relaxation capability makes a significant difference
in mesh quality improvement at the curved boundaries.}

\begin{figure}[bth]
\begin{center}
$\begin{array}{cc}
    \includegraphics[width=0.23\textwidth]{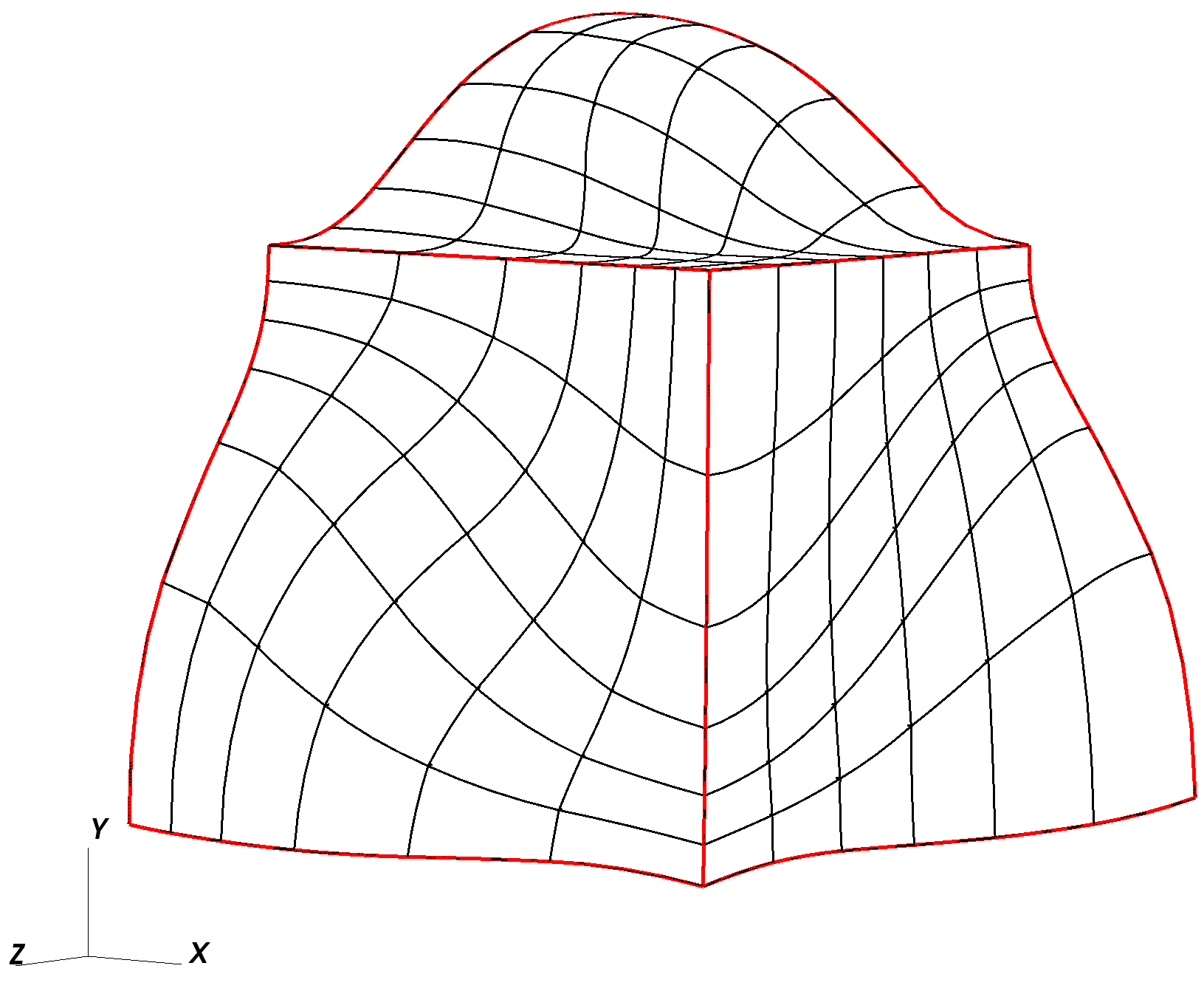} &
    \includegraphics[width=0.23\textwidth]{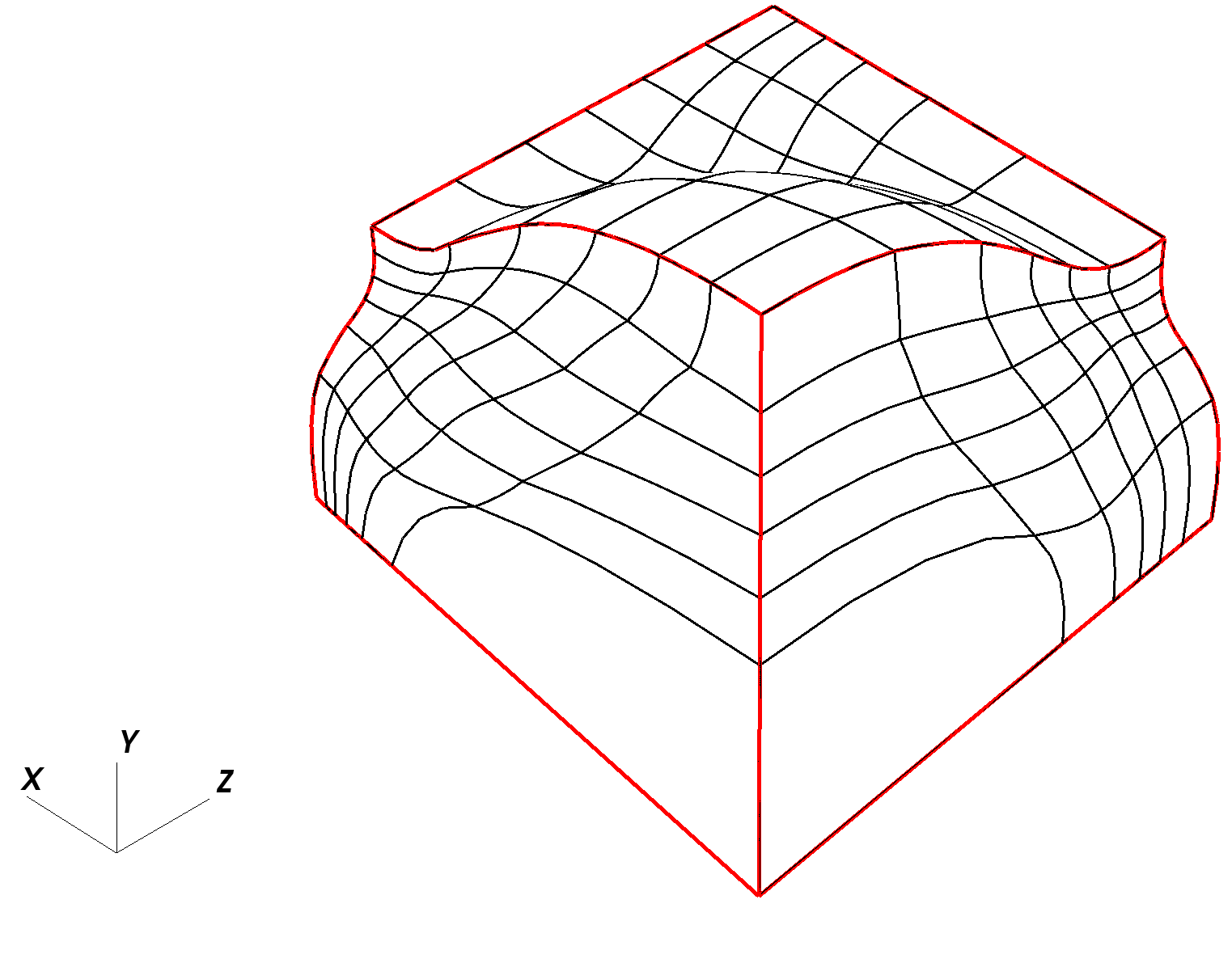} \\
\multicolumn{2}{c}{\textrm{(a) Initial mesh}} \\
    \includegraphics[width=0.23\textwidth]{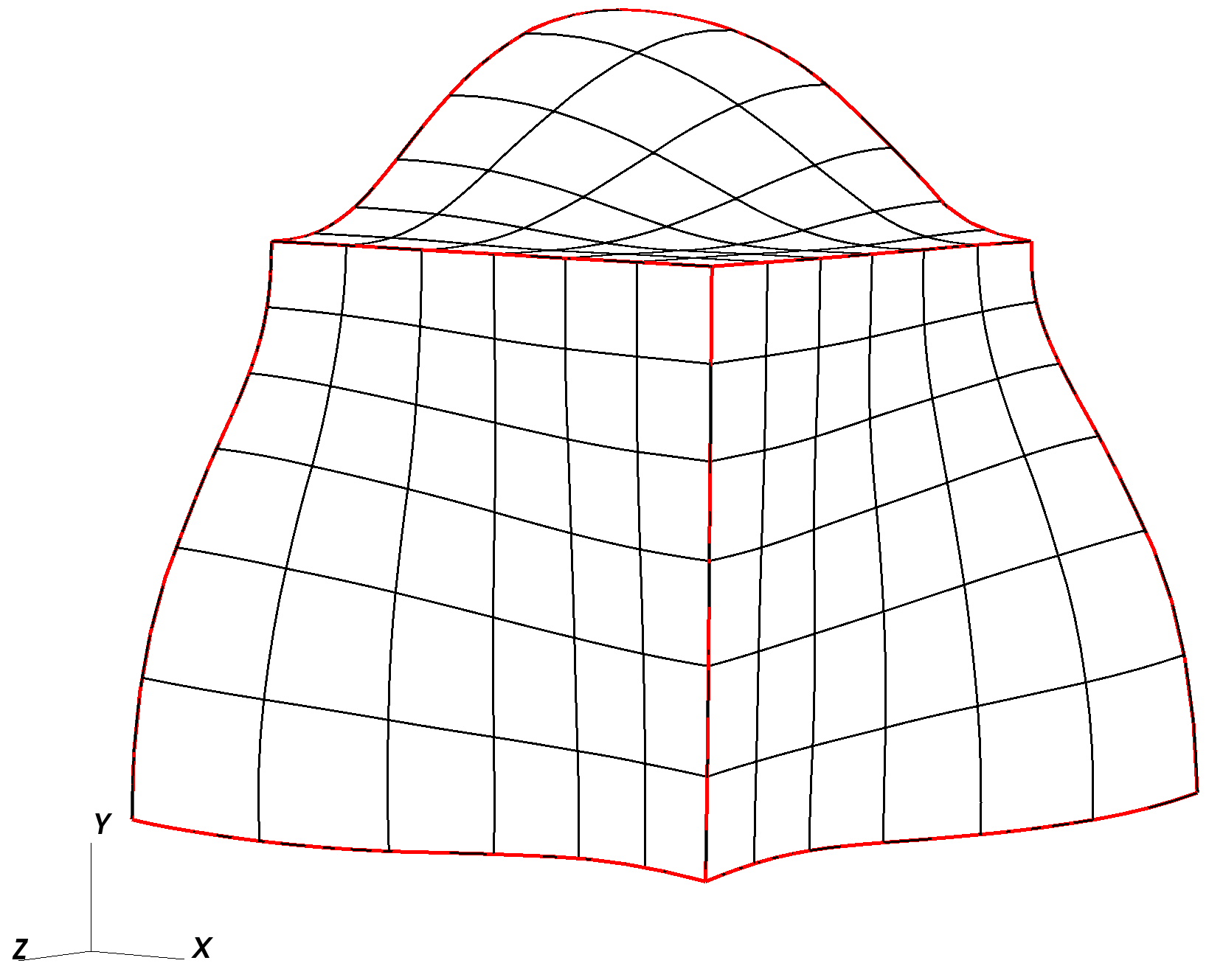} &
    \includegraphics[width=0.23\textwidth]{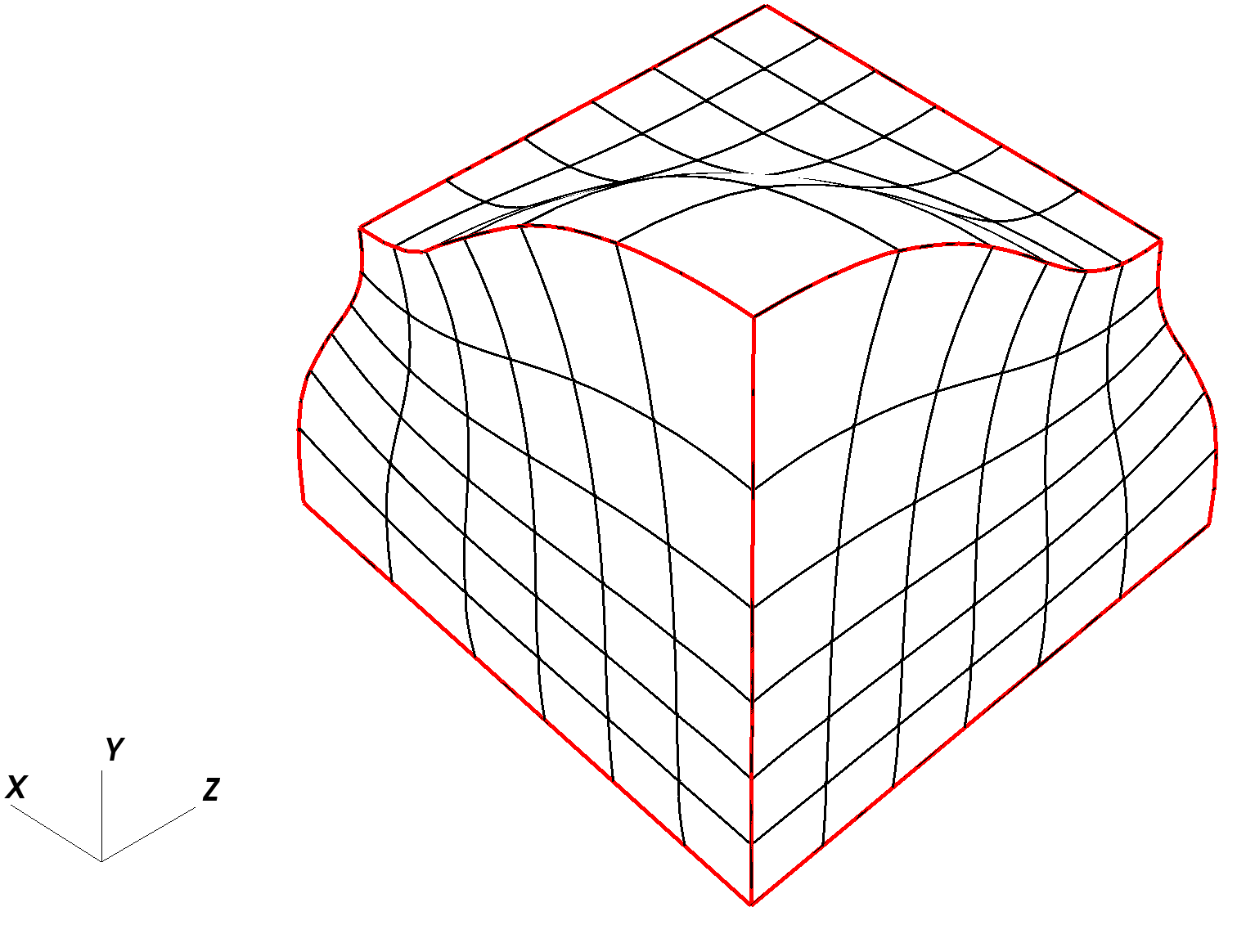} \\
\multicolumn{2}{c}{\textrm{(b) Optimized mesh}} \\
\end{array}$
\end{center}
\caption{$r$-adaptivity with geometry-preserving tangential relaxation for a 3D mesh.}
\label{fig:3D_ale}
\end{figure}

A similar 3D example with quadratic hexahedra is shown in Figure \ref{fig:3D_ale}.
The initial mesh is generated by performing a Lagrangian simulation of a shock wave,
initialized on a uniform Cartesian mesh. The shock originates at the corner
and deforms the domain while the boundaries are left free to move.
The deformed mesh is shown in Figure \ref{fig:3D_ale}(a). The lower bound on its Jacobian determinant is $6 \times 10^{-4}$, and the TMOP functional $\int \mu_{301}(T)$ is 225.4. Here, $\mu_{301}(T) = \frac{||T||_F ||T^{-1}||_F}{3} - 1$ is a shape metric. Geometric integrity of the domain boundaries is preserved during tangential relaxation by identifying mesh nodes that are on geometric edges and surfaces, as described in Section \ref{subsec:tangentialrelaxation}. These geometric edges are highlighted in red. Figure \ref{fig:3D_ale}(b) shows the optimized mesh with the lower bound on the Jacobian determinant increasing to $\ua (\bx) = 2.2 \times 10^{-3}$ and the TMOP functional decreasing to 3.34. This corresponds to a significant improvement in mesh quality, as evident in Figure \ref{fig:3D_ale}(b), where the mesh elements are more uniformly distributed and the aspect ratio is improved.

\section{Summary and Future Work}
\label{sec:conc}
We have extended the Target-Matrix Optimization Paradigm by developing a robust
framework for high-order mesh $r$-adaptivity that integrates tangential relaxation with provable guarantees on mesh validity, thereby improving mesh quality while ensuring positivity everywhere in the domain.
The effectiveness of this approach is demonstrated on challenging 2D and 3D examples, including ALE-type deformations where geometric fidelity and robustness are critical.

Future work will focus on extending the methodology to more complex and realistic 3D meshes that arise in practical applications.
In particular, we will explore its integration into multi-material ALE simulations, where accurate mesh adaptation plays a central role in capturing evolving material interfaces and maintaining numerical stability under strong deformations.

\section*{Acknowledgments}
This manuscript has been authored by Lawrence Livermore National Security, LLC under Contract No. DE-AC52-07NA27344 with the US. Department of Energy. The United States Government retains, and the publisher, by accepting the article for publication, acknowledges that the United States Government retains a non-exclusive, paid-up, irrevocable, world-wide license to publish or reproduce the published form of this manuscript, or allow others to do so, for United States Government purposes. LLNL-CONF-2011480.

\clearpage

\bibliographystyle{ieeetr}
\bibliography{lit}

\end{document}